\documentclass{elsart}

\usepackage{graphics,hyperref}
\usepackage{amsmath} 
\usepackage{amssymb}  
\usepackage{mathptmx}

\newcommand{\EEE}{{\mathcal E }}
\newcommand{\HHH}{{\mathcal H }}
\newcommand{\LLL}{{\mathcal L }}

\newcommand{\VVV}{{\mathcal  V}}
\newcommand{\PP}{{\mathbb P}}

\newcommand{\RR}{{\mathbb R}}
\newcommand{\NN}{{\mathbb N}}

\newcommand{\SSS}{{\mathbb  S}}
\newcommand{\CC}{{\mathbb C}}

\newcommand{\dotex}{{\frac{d}{dt}}}

\newcommand{\her}[2]{\left\langle#1~|~#2\right\rangle}

\begin{document}

\begin{frontmatter}


\title{Lyapunov control of a quantum particle in a decaying potential}
\author{Mazyar Mirrahimi}                     
\ead{mazyar.mirrahimi@inria.fr}

\thanks{This work was partly supported by the ``Agence Nationale de
la Recherche'' (ANR), Projet Blanc CQUID number 06-3-13957.}
\address{INRIA Rocquencourt, B.P. 105, Domaine de Voluceau, 78153
Le Chesnay Cedex, France.}


\begin{abstract}
A Lyapunov-based approach for the trajectory generation of an
$N$-dimensional Schr{\"o}dinger equation in whole $\RR^N$ is
proposed. For the case of a quantum particle in an $N$-dimensional
decaying potential the convergence is precisely analyzed. The free
system admitting a mixed spectrum, the dispersion through the
absolutely continuous part is the main obstacle to ensure such a
stabilization result. Whenever, the system is completely initialized
in the discrete part of the spectrum, a Lyapunov strategy encoding
both the distance with respect to the target state and the
penalization of the passage through the continuous part of the
spectrum, ensures the approximate stabilization.
\end{abstract}

\end{frontmatter}

\section{Introduction}\label{sec:intro}
\subsection{Main results}\label{ssec:intromain}
We consider a quantum particle in an $N$-dimensional space, with  a
potential $V(x)$, and coupled to an external (laser) field $t\mapsto
u(t)\in\RR$ through its dipole moment $\mu(x)$. Under appropriate
change of scales, the system's wavefunction evolves following the
Schr\"odinger equation
\begin{alignat}{2}
i\frac{\partial \Psi}{\partial t}(t,x) &=
-\triangle\Psi(t,x)+(V(x)-u(t)\mu(x))\Psi(t,x),\qquad
x\in\RR^N,\label{eq:main}\\
\Psi(0,x) &= \Psi_0(x).\label{eq:initial}
\end{alignat}
This is a bilinear control system, denoted by ($\Sigma$), where
\begin{itemize}
  \item the control is the external field $u:\RR_+\rightarrow\RR$,
  \item the state is the wave function
  $\Psi:\RR_+\times\RR^N\rightarrow \CC$ with $\Psi(t)\in\SSS$ for
  every $t\geq 0$, with $\SSS=\{\varphi\in
  L^2(\RR^N;\CC)~|~\|\varphi\|_{L^2}=1\}$.
\end{itemize}
We distinguish between four different situations: dimension $N\geq
4$; dimension $N=3$; dimension $N=2$; and dimension $N=1$. For each
of these cases, we will assume some appropriate decay assumptions
for the potential $V(x)$. Indeed, through this paper, we will assume
the following assumption:
\begin{description}
  \item[Decay assumption \textbf{(A)}] We assume for the potential $V$ that zero is neither an eigenvalue
nor a resonance of the Hamiltonian $H_0=-\triangle+V$. Furthermore,
we assume one of the following assumptions (depending on the space
dimension $N$)
\begin{itemize}
  \item $N=1$: $(1+|x|)V\in L^1(\RR)$~\cite{goldberg-schlag-04-1};
  \item $N=2$: $|V(x)|\leq C(1+|x|)^{-3-\epsilon}$~\cite{schlag-05};
  \item $N=3$: $V\in L^{\frac{3}{2}-\epsilon}(\RR^3)\cap
  L^{\frac{3}{2}+\epsilon}(\RR^3)$~\cite{goldberg-06};
  \item $N\geq4$: $\widehat V\in L^1$ and $(1+|x|^2)^{\gamma/2}V(x)$
  is a bounded operator on the sobolev space $H^\nu$ for some
  $\nu>0$ and $\gamma>n+4$~\cite{journe-soffer-sogge-91}.
\end{itemize}
\end{description}

A brief discussion on the origin of the above assumption on the
potential $V$ is provided in the Subsection~\ref{ssec:decay}. As one
will see these decay assumptions are chosen to assure relevant
dispersive estimates.

Furthermore, note that, under the decay assumption \textbf{(A)} on
the potential $V$, the free Hamiltonian $H_0=-\triangle+V(x)$ admits
a mixed spectrum:
$$
\sigma(H_0)=\sigma_{\text{disc}}(H_0)\cap \sigma_{ac}(H_0),
$$
where the discrete spectrum $\sigma_{\text{disc}}(H_0)$ contains a
finite number of eigenvalues of finite multiplicities and the
essential spectrum is actually an absolutely continuous spectrum
$\sigma_{ac}(H_0)=[0,\infty)$. Under the decay assumption
\textbf{(A)}, this decomposition of the spectrum for the 1D case is
a classical result of the earliest days of quantum mechanics (in
fact one only needs $V\in L^1(\RR)$, see e.g.~\cite{reed-simon4},
Sec. XIII.4). For the 2D case, one can find a proof
in~\cite{stoiciu-04}. The 3D case has been proven
in~\cite{goldberg-schlag-04-2}. Finally, the decomposition for the
$N$-dimensional case, with $N>3$, is a classical result as the
potential is a short range potential in the sense of
Agmon~\cite{Agmon-75}.

Concerning the bound states, $\{\phi_j\}_{j=0}^M$, we know that
$\phi_j\in H^2(\RR^N,\CC)$. Moreover, the decay assumption
\textbf{(A)} on the potential implies that $V\in L^1_{\text{loc}}$
and $V_-\in M_{\text{loc}}$ the local Stummel class
(see~\cite{agmon-book}, page 8, for a definition). This ensures the
exponential decay of the eigenfunctions $\{\phi_j\}_{j=0}^M$ (see
e.g.~\cite{agmon-book}, page 55, Corollary 4.2).

Let us recall the following classical existence and uniqueness
result for the open-loop system~\eqref{eq:main}-\eqref{eq:initial}.
A proof of this result is given in the Appendix.
\begin{prop}\label{prop:well-posed}
Let the potential $V(x)$ satisfy the decay assumption \textbf{(A)}
and consider $\mu\in L^\infty(\RR^N,\RR)$. Let $\Psi_0\in\SSS$,
$T>0$ and $u\in \CC^0([0,T],\RR)$. There exists a unique weak
solution of~\eqref{eq:main}-\eqref{eq:initial}, i.e. a function
$\Psi\in C^0([0,T],\SSS)\cap C^1([0,T],H^{-2}(\RR^N,\CC))$ such that
\begin{multline}\label{eq:mild}
\Psi(t)=e^{-iH_0 t}\Psi(0)\\+i \int_0^t e^{-iH_0 (t-s)} u(s) \mu(x)
\Psi(s)ds\qquad \text{ in } L^2(\RR^N,\CC) \text{ for } t\in[0,T],
\end{multline}
and then~\eqref{eq:main} holds in $H^{-2}(\RR^N,\CC)$.

If, moreover, $\Psi_0\in H^2(\RR^N,\CC)$ and multiplication by
$\mu(x)$ defines a bounded operator over $H^2(\RR^N,\RR)$, then
$\Psi$ is a strong solution, i.e. $\Psi\in
C^0([0,T],H^2(\RR^N,\CC))$ $\cap C^1([0,T],L^2(\RR^N,\CC))$, the
equation~\eqref{eq:main} holds in $L^2(\RR^N,\CC)$ for $t\in[0,T]$
and the initial condition~\eqref{eq:initial} holds in
$H^2(\RR^N,\CC)$.

The weak (resp. strong) solution is continuous with respect to the
initial condition for the $C^0([0,T],L^2)$-topology (resp.
$C^0([0,T],H^2)$-topology).
\end{prop}

Assuming the potential $V(x)$ such that the discrete spectrum
$\sigma_{\text{disc}}(H_0)$ is non-empty, we are interested here in
stabilizing one of the eigenfunctions in this discrete part. Fixing
$\epsilon>0$ to be a small positive constant and considering $\phi$
to be a normalized eigenfunction in this discrete part, we are
interested in designing a feedback law $u_{\epsilon}(\Psi)$ such
that, the solution $\Psi(t,x)$ of~\eqref{eq:main}-\eqref{eq:initial}
satisfies
\begin{equation}\label{eq:liminf}
\liminf_{t\rightarrow\infty}
|\her{\Psi(t,x)}{\phi(x)}|^2>1-\epsilon.
\end{equation}
Here
$$
\her{\xi}{\zeta}=\int_{\RR^N} \xi(x)\overline\zeta(x) dx,
$$
denotes the Hermitian product of $L^2(\RR^N,\CC)$.  Note that,
$\Psi$ and $\phi$ living on the unit sphere $\SSS$ of
$L^2(\RR^N,\CC)$, the limit~\eqref{eq:liminf} denotes the
approximate stabilization of the eigenfunction $\phi(x)$.

Note that, even though the feedback stabilization of a quantum
system necessitates more complicated models taking into account the
measurement backaction on the system (see
e.g.~\cite{haroche-CDF,vanhandel-et-al-05,mirrahimi-vanhandel-05}),
the kind of strategy considered in this paper can be helpful for the
open-loop control of closed quantum systems. Indeed, one can apply
the stabilization techniques for the Schr\"odinger equation in
simulation and retrieve the control signal that will be then applied
in open-loop on the real physical system. As it will be detailed
below, in the bibliographic overview, such kind of strategy has been
widely used in the context of finite dimensional quantum systems.

The main result of this article is the following one.
\begin{thm}\label{thm:gas}
Consider the Schr{\"o}dinger
equation~\eqref{eq:main}-~\eqref{eq:initial}. We suppose the
potential $V(x)$ to satisfy the decay assumption \textbf{(A)} and we
take $\mu\in \LLL(\RR^N)\cap L^\infty(\RR^N)$. We assume the
discrete spectrum $\sigma_{\text{disc}}$ of $H_0=-\triangle+V(x)$ to
be non-empty. We consider moreover the following assumptions:
\begin{description}
    \item[A1] $\Psi_0=\sum_{j=0}^M
    \alpha_j\phi_j$ where $\{\phi_j\}_{j=0}^M$ are different
    normalized eigenfunctions in the discrete spectrum of $H_0$.
    \item[A2] the coefficient $\alpha_0$ corresponding
    to the population of the eigenfunction $\phi_0$ in the initial
    condition $\Psi_0$ is non-zero: $\alpha_0\neq 0$.
    \item[A3] the Hamiltonian $H_0$ admits
    non-degenerate transitions: $\lambda_{j_1}-\lambda_{k_1}\neq
    \lambda_{j_2}-\lambda_{k_2}$ for $(j_1,k_1)\neq (j_2,k_2)$ and
    where $\{\lambda_j\}_{j=0}^M$ are different eigenvalues of the
    Hamiltonian $H_0$;
    \item[A4] the interaction Hamiltonian $\mu(x)$ ensures simple
    transitions between all eigenfunctions of $H_0$:
    $$
    \her{\mu\phi_j}{\phi_k}\neq 0 \qquad \forall j\neq k
    \in\{0,1,...,M\}.
    $$
\end{description}
Then for any $\epsilon >0$, there exists a feedback law
$u(t)=u_\epsilon(\Psi(t))$ (that we will construct explicitly), such
that the closed-loop system admits a unique weak solution in
$C^0([0,T],\SSS)\cap C^1([0,T],H^{-2}(\RR^N,\CC))$. Moreover the
state of the system ends up reaching a population more than
$(1-\epsilon)$ in the eigenfunction $\phi_0$ (approximate
stabilization):
$$
\liminf_{t\rightarrow
\infty}|\her{\Psi(t,x)}{\phi_0(x)}|^2>1-\epsilon.
$$
If, moreover multiplication by $\mu(x)$ defines a bounded operator
over $H^2(\RR^N)$, then $\Psi$ is a strong solution, i.e. $\Psi\in
C^0([0,T],H^2$ $(\RR^N,\CC))\cap C^1([0,T],L^2(\RR^N,\CC))$.
\end{thm}
\begin{rem}\label{rem:Lp-}
In this Theorem $\LLL$ denotes $\bigcup_{p\geq2} L^p(\RR^N)$.
\end{rem}
\begin{rem}\label{eq:initial2}
Note that, as the initial state is a linear combination of the bound
states, we have in particular $\Psi_0\in H^2$ and decays
exponentially.
\end{rem}
\begin{rem}\label{rem:rev}
Note that, here a finite dimensional approximation of the system by
removing the continuous part of the spectrum is not sufficient to
treat the stabilization problem. In fact, even if the system is
initialized in the discrete part of the spectrum (as assumed in
\textbf{A1}), the interaction Hamiltonian $\mu$ will make the
solution leave this discrete part. The state of the system will
therefore leave this subspace just after the initial time.
\end{rem}
The assumptions \textbf{A1} through \textbf{A4} can be relaxed
significantly. However, as the final result with the relaxed
assumptions may seem too complicated, we will discuss this
relaxations, separately, in Section~\ref{sec:relax}.

\subsection{A brief bibliography}\label{ssec:biblio}
The controllability of a finite dimensional quantum system,
$i\frac{d}{dt}\Psi=(H_0+u(t)~H_1)\Psi$ where $\Psi\in \CC^N$ and
$H_0$ and $H_1$ are $N\times N$ Hermitian matrices with coefficients
in $\CC$ has been completely
explored~\cite{sussmann-jurdjevic-72,ramakrishna-et-al-95,albertini-et-al-03,altafini-JMP-02,turinici-rabitz-03}.
However, this does not guarantee the simplicity of the trajectory
generation. Very often the chemists formulate the task of the
open-loop control as a cost functional to be minimized. Optimal
control techniques (see e.g.,~\cite{hbref17}) and iterative
stochastic techniques (e.g, genetic
algorithms~\cite{turinici-control1}) are then two classes of
approaches which are most commonly used for this task.

When some non-degeneracy assumptions concerning the linearized
system are satisfied,~\cite{mirrahimi-et-al2-04} provides another
method based on Lyapunov techniques for generating trajectories. The
relevance of such a method for the control of chemical models has
been studied in~\cite{mirrahimi-et-al04}. Since measurement and
feedback in quantum systems lead to much more complicated models and
dynamics than the Schr\"{o}dinger
equation~\cite{haroche-CDF,mirrahimi-vanhandel-05}, the
stabilization techniques presented in~\cite{mirrahimi-et-al2-04} are
only used for generating open-loop control laws. Simulating the
closed-loop system, we obtain a control signal which can be used in
open-loop for the physical system. Such kind of strategy has already
been applied widely in this
framework~\cite{rabitz-tracking-95,sugawara-JCP03}.

The situation is much more difficult when we consider an infinite
dimensional configuration. Concerning the controllability problem,
very few results are
available~\cite{turinici-cdc00,beauchard-04,beauchard-coron-05}.
In~\cite{beauchard-04,beauchard-coron-05} the controllability of a
particle in a moving one dimensional quantum box has been studied. A
local controllability result is therefore provided using the return
method~\cite{coron-mcss-92}. In~\cite{chambrion-et-al-08}, applying
some geometric control tools, the authors provide a quite general
result concerning the controllability of discrete-spectrum
Schr\"odinger equation. Finally, in~\cite{tarn-et-al-00}, the
authors consider the controllability of some particular
Schr\"odinger equations with continuous spectra.

Concerning the trajectory generation problem for infinite
dimensional systems still much less results are available. The very
few existing literature is mostly based on the use of the optimal
control techniques~\cite{baudouin-et-al-05,baudouin-salomon-06}. The
simplicity of the feedback law found by the Lyapunov techniques
in~\cite{mirrahimi-et-al2-04,beauchard-et-al07} suggests the use of
the same approach for infinite dimensional configurations. However,
an extension of the convergence analysis to the PDE configuration is
not at all a trivial problem. Indeed, it requires the
pre-compactness of the closed-loop trajectories, a property that is
difficult to prove in infinite dimension. This strategy is used, for
example in \cite{JMC-BAN}.

Let us mention some strategies for proving the stabilization of
infinite dimensional control systems. One can try to build a
feedback law for which one has a strict Lyapunov function. This
strategy is used, for example, for hyperbolic systems of
conservation laws in \cite{JMC-BAN-Lyapunov-strict}, for the 2-D
incompressible Euler equation in a simply connected domain in
\cite{coron-siam-99}, see also \cite{Glass-Stab-Euler2D} for the
multi-connected case. For systems having a non controllable
linearized system around the equilibrium considered, the return
method often provides good results, see for example
\cite{coron-mcss-92} for controllable systems without drift and
\cite{Glass-Camassa}) for Camassa-Holm equation. In the end, we
refer to \cite{JMC-book} for a pedagogical presentation of
strategies for the proof of stabilization of PDE control systems.

In this paper, we propose a Lyapunov-based method to approximately
stabilize a particle in an $N$ dimensional decaying potential under
some relevant assumptions. We assume that the system is initialized
in the finite dimensional discrete part of the spectrum. Then, the
idea consists in proposing a Lyapunov function which encodes both
the distance with respect to the target state and the necessity of
remaining in the discrete part of the spectrum. In this way, we
prevent the possibility of the ``mass lost phenomenon'' at infinity.
Finally, applying some dispersive estimates of Strichartz type, we
ensure the approximate stabilization of an arbitrary eigenfunction
in the discrete part of the spectrum.

The ideas of this paper (a short and simplified version is already
published as a communication~\cite{mirrahimi-cdc06}) have been
recently adapted to the case of a quantum particle in an infinite
potential well~\cite{beauchard-mirrahimi-07}.
In~\cite{beauchard-mirrahimi-07}, as we are dealing with a pure
discrete spectrum, much less restrictive assumptions are needed to
ensure the approximate stabilization of the system.

As it can be remarked through the bibliography, except for a very
few results~\cite{tarn-et-al-00}, all the previous work on the
control of the infinite dimensional quantum systems deal with
discrete-spectrum Schr\"{o}dinger equations. It seems that the
techniques of this paper and the possibility of the relaxations,
explained in Section~\ref{sec:relax}, can open a new gateway to
investigate this class of quantum systems.
\subsection{Free dynamics and dispersive estimates}\label{ssec:decay}
Before treating the control problem, let us have a look at the
behavior of the system in the absence of the control field
($u(t)=0$). We will denote by $S(t)=\exp(-itH_0)$ the
$C_0$-semigroup on $L^2(\RR^N,\CC)$ spanned by the infinitesimal
generator $(-\triangle+V(x))/i$. Note in particular that, $S(t)$
induces an isometry over $L^2(\RR^N,\CC)$:
$\|S(t)\psi\|_{L^2}=\|\psi\|_{L^2}$.

Moreover, we denote by $\PP_{\text{disc}}$ the projection operator
over the discrete subspace generated by the bound states and defined
on $L^2(\RR^N,\CC)$. Finally, $\PP_{\text{ac}}$ denotes the
projection over the orthogonal subspace:
$\PP_{\text{ac}}=Id-\PP_{\text{disc}}$.

The discrete part of the freely evolving solution $\PP_{\text{disc}}
S(t)\Psi_0$ represents a quasi-periodic behavior:
$$
\Psi_{0,\text{disc}}=\PP_{\text{disc}}\Psi_0=\sum_{j=0}^M \alpha_j
\phi_j(x)\qquad \Rightarrow\qquad \PP_{\text{disc}} S(t) \Psi_0=
\sum_{j=0}^M \alpha_j e^{-i\lambda_j t} \phi_j(x).
$$
The continuous part, however, represents a dispersive behavior. In
this subsection, we provide a very brief overview of the dispersive
estimates and in particular the ones we use in this paper.

In heuristic terms, for the potential-free problem $V\equiv 0$, the
explicit solution
$$
(e^{it\triangle}\psi)(x)=C_N~t^{-N/2}\int_{\RR^N}
e^{i\frac{|x-y|^2}{4t}}\psi(y)dy,
$$
implies the dispersive estimate~\cite{strichartz}
$$
\sup_{t>0} |t|^{N/2}\|e^{it\triangle}\psi\|_{L^\infty(\RR^N)}\leq
\|\psi\|_{L^1(\RR^N)}\quad \forall \psi\in L^1(\RR^N)\cap
L^2(\RR^N).
$$
For general $V\neq 0$, no explicit solutions are available and
therefore one needs to proceed differently. Consider the perturbed
Hamiltonian $H_0= -\triangle + V$, we seek to prove similar
estimates on the time evolution operator
$S(t)\PP_{\text{ac}}=e^{-itH_0}\PP_{\text{ac}}$. The projection onto
the absolutely continuous spectrum of $H_0$ is needed to eliminate
bound states which do not decay over any length of time. We,
therefore, have the following dispersive estimate:
\begin{thm}\label{thm:disp}
Under the decay assumption \textbf{(A)} on the potential $V$, we
have
\begin{equation}\label{eq:disp}
\|S(t)\PP_{\text{ac}}\|_{1\rightarrow\infty}\leq |t|^{-\frac{N}{2}}.
\end{equation}
\end{thm}
Such dispersive estimates have a long history. For exponentially
decaying potentials, Rauch~\cite{rauch-78} proved dispersive bounds
in exponentially weighted $L^2$-spaces. Jensen and
Kato~\cite{jensen-kato} replaced exponential with polynomial decay
and obtained asymptotic expansions of $e^{-itH_0}$ (in terms of
powers of $t$) in the usual weighted $L^{2,\sigma}$ spaces. The
first authors to address a dispersive estimate of the
form~\eqref{eq:disp} were Journ\'ee, Soffer, and
Sogge~\cite{journe-soffer-sogge-91}. They were able to prove the
dispersive estimate~\eqref{eq:disp} under the fourth case of the
decay assumption \textbf{(A)} for $N\geq 3$.

Concerning the case $N=3$, following a large amount of
results~\cite{yajima-95,rodnianski-schlag-04,goldberg-schlag-04-1},
finally Goldberg~\cite{goldberg-06} proved the dispersive
estimate~\eqref{eq:disp} under the third case of the decay
assumption \textbf{(A)}. In contrast, trying to adapt these results
to higher dimensions has lead Goldberg and
Visan~\cite{goldberg-visan-06} to show that for $N \geq
4$,~\eqref{eq:disp} fails unless $V$ has some amount of regularity,
i.e., decay alone is insufficient for~\eqref{eq:disp} to hold if
$N\geq 4$.

The one-dimensional case was open until recently.
Weder~\cite{weder-00} proved a version of Theorem~\ref{thm:disp}
under the stronger assumption that $\int_{-\infty}^\infty
|V(x)|(1+|x|)^{3/2+\epsilon}dx<\infty$. Finally, in a similar way
to~\cite{weder-00}, Goldberg and Schlag~\cite{goldberg-schlag-04-1}
were able to prove~\eqref{eq:disp} under the first case of the decay
assumption \textbf{(A)}.

Finally, concerning the two-dimensional case,
Yajima~\cite{yajima-99} and Jensen, Yajima~\cite{jensen-yajima}
proved the $L^p(\RR^2)$ bounded-ness of the wave operators under
stronger decay assumptions on $V(x)$ (than the decay assumption
\textbf{(A)}), but only for $1 < p < \infty$. Hence their result
does not imply~\eqref{eq:disp}, but $L^p\rightarrow L^{p'}$
estimates for $1<p\leq 2$. The first paper to provide an
$L^1\rightarrow L^\infty$ dispersive estimate of the
form~\eqref{eq:disp} in two dimensions was that of
Schlag~\cite{schlag-05} that proves~\eqref{eq:disp} under the second
case of the decay assumption \textbf{(A)}.

Note that, interpolating with the $L^2$-bound
$\|e^{-itH_0}\PP_{\text{ac}}\psi\|_{L^2}\leq \|\psi\|_{L^2}$, we
have
\begin{cor}\label{cor:Lpq}
Under the decay assumption \textbf{(A)} on the potential $V$, we
have
\begin{equation}\label{eq:disp2}
\sup_{t>0}
|t|^{N\left(\frac{1}{2}-\frac{1}{p}\right)}\|S(t)\PP_{\text{ac}}\psi\|_{L^{p'}}\leq
\|\psi\|_{L^p}\quad\text{ for all }\psi\in L^1(\RR^N)\cap
L^2(\RR^N),
\end{equation}
where $1\leq p\leq 2$ and $\frac{1}{p}+\frac{1}{p'}=1$.
\end{cor}
Furthermore, through a $T^*T$ argument,~\eqref{eq:disp2} leads to
the class of Strichartz estimates (see e.g.~\cite{keel-tao}):
\begin{thm}\label{thm:strichartz}
Under the decay assumption \textbf{(A)} on the potential $V$, we
have
\begin{equation}\label{eq:disp3}
\|S(t)\PP_{\text{ac}}\psi\|_{L^q_t(L^p_x)}\leq C\|\psi\|_{L^2},
\text{ for all } \frac{2}{q}+\frac{N}{p}=\frac{N}{2},\quad 2<q\leq
\infty.
\end{equation}
\end{thm}
\subsection{Structure of the paper}\label{ssec:struct}
The rest of the paper is organized as follows. In
Section~\ref{sec:heur}, we will provide the heuristic of the proof
of the Theorem~\ref{thm:gas}. In this aim, we will first announce a
new Theorem~\ref{thm:ggg} providing the same result as the
Theorem~\ref{thm:gas} but under some more restrictive assumptions.
After discussing heuristically the proof of this new theorem, we
will give the elements to extend the proof to that of the
Theorem~\ref{thm:gas}.

The Theorem~\ref{thm:ggg}, will be proved in
Section~\ref{sec:proof1}. Through a rather simple change in the
feedback law, we will be able to extend this proof to that of the
Theorem~\ref{thm:gas}. This will be addressed in
Section~\ref{sec:proof2}.

As it can be seen, the assumptions of the Theorem~\ref{thm:gas} may
still seem too restrictive. However, through some arguments based on
the analytic perturbation of linear operators and the quantum
adiabatic theory, we are able to relax significantly these
assumptions. This will be treated in Section~\ref{sec:relax}.
\section{Heuristic of the proof}\label{sec:heur}
From now on, we will assume that the system is initially prepared in
a purely discrete state:
$$
\Psi_0=\Psi_{0,\text{disc}}\in\EEE_{\text{disc}},
$$
where $\EEE_{\text{disc}}$ (resp. $\EEE_{\text{ac}}$) denotes
$\text{Range}(\PP_{\text{disc}})$ (resp.
$\text{Range}(\PP_{\text{ac}})$). The control task is to steer the
systems state in the eigenspace corresponding to an eigenfunction
$\phi_0$ of the free Hamiltonian. Note that this eigenfunction
$\phi_0$ can be any eigenfunction in the discrete part of the
spectrum and does not have to be the ground state. During the
control process the system might and will cross the continuum
$\EEE_{\text{ac}}$.

Following the stabilization results for the finite dimensional
systems \cite{mirrahimi-et-al2-04,beauchard-et-al07}, a first
approach for this control problem might be to consider the simple
Lyapunov function
$$\widetilde\VVV(\Psi)=1-|\her{\Psi}{\phi_0}|^2.$$
The fact that $\Psi$ and $\phi_0$ are both normalized, together with
the Cauchy-Schwartz inequality, ensures that $\widetilde
\VVV(\Psi)\geq 0$. The feedback law will be given
by~\cite{beauchard-et-al07}:
\begin{equation}\label{eq:try}
\tilde u(\Psi)=\Im(\her{\mu\Psi}{\phi_0}\her{\phi_0}{\Psi}),
\end{equation}
where $\Im$ denotes the imaginary part of a complex number. A deep
analysis based on LaSalle type arguments shows that with such a
feedback strategy, one can not avoid phenomenons like mass lost at
infinity. The population of the state $\phi_0$ will surely keep
increasing during the evolution. But, in order to be able to apply
the LaSalle invariance principle for such infinite dimensional
system, one needs to ensure the pre-compactness of the trajectories
in $L^2(\RR^N,\CC)$. In the particular case, of the Schr\"odinger
equation with the decaying potential, considered in this paper and
with the feedback law~\eqref{eq:try}, one can not even hope to have
such a pre-compactness result. Indeed, as it has been said before,
while the population of the state $\phi_0$ keeps increasing through
the application of the feedback law~\eqref{eq:try}, during this same
period some of the population might go through the continuous part
of the spectrum. This population has then the possibility to
disperse rapidly (cf. Subsection~\ref{ssec:decay}) and so we might
have some un-controlled part of the $L^2$-norm which will be lost at
infinity.

The approach of this paper consists in avoiding the population to go
through the continuum while stabilizing the state around the target
state $\phi_0$. So, we consider  a Lyapunov function
$\VVV_{\epsilon}(\Psi)$ which encodes these both tasks:
\begin{equation}\label{eq:lyap}
\VVV_{\epsilon}(\Psi):=1-(1-\epsilon)\sum_{j=0}^M
|\her{\Psi}{\phi_j}|^2-\epsilon|\her{\Psi}{\phi_0}|^2,
\end{equation}
where $0<\epsilon\ll 1$ is a small positive constant. Such a
Lyapunov function clearly verifies:
\begin{equation}\label{eq:pos}
0\leq \VVV_\epsilon(\Psi) \qquad \text{and} \qquad
\VVV_\epsilon(\Psi)=0 \Leftrightarrow |\her{\Psi}{\phi_0}|=1.
\end{equation}
Here still, we have used the fact that $\Psi$ and $\phi_j$'s are all
normalized in $L^2(\RR^N,\CC)$. Moreover, as the system is initially
prepared in the discrete part of the spectrum, and as
$|\her{\Psi_0}{\phi_0}|>0$,
\begin{equation}\label{eq:initlyap}
\VVV_\epsilon(\Psi_0)=1-(1-\epsilon)-\epsilon|\her{\Psi_0}{\phi_0}|<
\epsilon.
\end{equation}
This Lyapunov function clearly encodes two tasks: 1- it prevents the
$L^2$-mass lost through the dispersion of the absolutely continuous
population; 2- it privileges the increase of the population in the
eigenfunction $\phi_0$.

By a simple computation we have,
\begin{equation}\label{eq:feedlyap}
\frac{d}{dt}\VVV_\epsilon(\Psi)=-u(t)[(1-\epsilon)\sum_{j=0}^M
\Im(\her{\mu\Psi}{\phi_j}\her{\phi_j}{\Psi})  +\epsilon
\Im(\her{\mu\Psi}{\phi_0}\her{\phi_0}{\Psi})].
\end{equation}
A natural choice is therefore to consider the feedback law:
\begin{equation}\label{eq:feed}
u(\Psi)=u_\epsilon(\Psi):=
c[(1-\epsilon)\sum_{j=0}^M\Im(\her{\mu\Psi}{\phi_j}\her{\phi_j}{\Psi})
+\epsilon\Im(\her{\mu\Psi}{\phi_0}\her{\phi_0}{\Psi})],
\end{equation}
where $c>0$ is a positive constant. Such a feedback law clearly
ensures the decrease of the Lyapunov function $\VVV_\epsilon$.
Looking at the structure of the Lyapunov function~\eqref{eq:lyap},
the feedback law~\eqref{eq:feed} penalizes strongly exiting from the
discrete part of the spectrum. Actually, as
$\VVV_\epsilon(\Psi_0)\leq \epsilon$, the decrease in this Lyapunov
function ensures that the population in the discrete part of the
spectrum will always remain more than $1-\epsilon$. Therefore in the
worst case, we will only have an $\epsilon$ $L^2$-norm which will be
lost by dispersing in the continuum.

At the same time, this feedback law~\eqref{eq:feed} slightly
encourages the increase in the population of the target state
$\phi_0$. It remains therefore to check whether this increase
actually provides some kind of convergence toward this eigenfunction
or not. This will be addressed in Section~\ref{sec:proof1}, where we
prove the following Theorem:
\begin{thm}\label{thm:ggg}
Consider the Schr{\"o}dinger
equation~\eqref{eq:main}-~\eqref{eq:initial}. Assume the space
dimension $N\geq 2$. We suppose the potential $V(x)$ to satisfy the
decay assumption \textbf{(A)} and we take $\mu\in L^{2N-}(\RR^N)\cap
L^\infty(\RR^N)$. We assume the discrete spectrum
$\sigma_{\text{disc}}$ of $H_0=-\triangle+V(x)$ to be non-empty. We
consider moreover the following assumptions:
\begin{description}
    \item[A1] $\Psi_0=\sum_{j=0}^M
    \alpha_j\phi_j$ where $\{\phi_j\}_{j=0}^M$ are different
    normalized eigenfunctions in the discrete spectrum of $H_0$.
    \item[A2] the coefficient $\alpha_0$ corresponding
    to the population of the eigenfunction $\phi_0$ in the initial
    condition $\Psi_0$ is non-zero: $\alpha_0\neq 0$.
    \item[A3] the Hamiltonian $H_0$ admits
    non-degenerate transitions: $\lambda_{j_1}-\lambda_{k_1}\neq
    \lambda_{j_2}-\lambda_{k_2}$ for $(j_1,k_1)\neq (j_2,k_2)$ and
    where $\{\lambda_j\}_{j=0}^M$ are different eigenvalues of the
    Hamiltonian $H_0$;
    \item[A4] the interaction Hamiltonian $\mu(x)$ ensures simple
    transitions between all eigenfunctions of $H_0$:
    $$
    \her{\mu\phi_j}{\phi_k}\neq 0 \qquad \forall j\neq k
    \in\{0,1,...,M\}.
    $$
\end{description}
Then for any $\epsilon >0$, applying the feedback law
$u(t)=u_\epsilon(\Psi(t))$ given by~\eqref{eq:feed}, the closed-loop
system admits a unique weak solution in $C^0([0,T],\SSS)\cap
C^1([0,T],H^{-2}(\RR^N,\CC))$. Moreover the state of the system ends
up reaching a population more than $(1-\epsilon)$ in the
eigenfunction $\phi_0$ (approximate stabilization):
$$
\liminf_{t\rightarrow
\infty}|\her{\Psi(t,x)}{\phi_0(x)}|^2>1-\epsilon.
$$
If, moreover, the multiplication by $\mu(x)$ defines a bounded
operator over $H^2(\RR^N)$, then $\Psi$ is a strong solution, i.e.
$\Psi\in C^0([0,T],H^2$ $(\RR^N,\CC))\cap
C^1([0,T],L^2(\RR^N,\CC))$.
\end{thm}
\begin{rem}\label{rem:ggg}
This theorem admits some more restrictive assumptions with respect
to the Theorem~\eqref{thm:gas}. In fact, we remove the 1D case and
we assume the interaction Hamiltonian $\mu$ to be in a smaller space
$L^{2N-}(\RR^N)\cap L^\infty(\RR^N)$, where
$$
L^{2N-}(\RR^N)=\bigcup_{2\leq p<2N} L^p,
$$
and therefore, $L^{2N-}(\RR^N)\subset \LLL(\RR^N)$.
\end{rem}
The Theorem~\ref{thm:ggg} will be proved by studying the $L^2$-weak
limit of $\Psi(t)$ for $t\rightarrow\infty$. Namely, let
$(t_n)_{n\in\NN}$ be an increasing sequence of positive real numbers
such that $t_n\rightarrow \infty$. Since $\|\Psi(t_n)\|_{L^2}=1$,
there exists $\Psi_\infty\in L^2(\RR^N,\CC)$ such that, up to a
subsequence, $\Psi(t_n)\rightharpoonup \Psi_\infty$ weakly in
$L^2(\RR^N,\CC)$. Furthermore, through the dispersive estimates of
the Subsection~\ref{ssec:decay}, we will provide a strong
convergence result with respect  to the semi-norm
$\|\psi\|_{\HHH}=\max(\|\PP_{\text{disc}}\psi\|_{L^2},\|\mu\psi\|_{L^2})$.
Through such a strong convergence and the Assumptions \textbf{A3}
and $\textbf{A4}$ of the Theorem, we will prove that
$\Psi_\infty=\beta\phi_0$, where $\beta\in\CC$ and $ |\beta|\leq 1$.
Through some further investigations and applying the
Assumption~\textbf{A1} we will be able to show that $|\beta|^2\geq
1-\epsilon$ and this will finish the proof of the
Theorem~\ref{thm:ggg}.

A deep study of the proof of the Theorem~\ref{thm:ggg}, shows that
the new restrictions (with respect to the Theorem~\ref{thm:gas}) may
be removed if we could ensure the belonging of the feedback law
$u_\epsilon(\Psi(t))$ to the space $L_t^{1+\delta}$ for $\delta \in
(0,1]$. In fact, the feedback law~\eqref{eq:feed} only belongs to
the space $L^2_t$ as
$\frac{d\VVV_\epsilon}{dt}=-\frac{1}{c}u_\epsilon^2$. However, we
may improve this through the following change of the feedback law:
\begin{equation}\label{eq:feed2}
u_{\epsilon,\alpha}(\Psi)=c f(\Psi)|f(\Psi)|^\alpha,
\end{equation}
where
$$
f(\Psi):=[(1-\epsilon)\sum_{j=0}^M\Im(\her{\mu\Psi}{\phi_j}\her{\phi_j}{\Psi})
+\epsilon\Im(\her{\mu\Psi}{\phi_0}\her{\phi_0}{\Psi})],
$$
and $\alpha\geq 0$ and $c>0$. This choice of the feedback law
implies
$$
\frac{d}{dt}\VVV_\epsilon=-\frac{c^{1+\alpha}}{c^{2+\alpha}}|u_{\epsilon,\alpha}|^{\frac{2+\alpha}{1+\alpha}},
$$
and therefore $u_{\epsilon,\alpha}\in
L_t^{\frac{2+\alpha}{1+\alpha}}$. As $\alpha\rightarrow\infty$ this
ensures that the feedback law $u_{\epsilon,\alpha}$ belongs to
$L_t^{1+\delta}$ for any $\delta\in(0,1]$.
\section{Proof of the Theorem~\ref{thm:ggg}}\label{sec:proof1}
We proceed the proof of the Theorem~\ref{thm:ggg} in 3 steps: 1- we
prove the well-posedness of the closed-loop system; 2- we prove the
existence of an asymptotic regime in some appropriate Hilbert space
and we characterize the weak $\omega$-limit set, i.e. the set of the
functions $\psi_\infty$ in $L^2(\RR^N)$ such that there exists a
sequence of times $\{t_n\}_{n=1}^\infty\nearrow\infty$ such that
$\Psi(t_n)\rightharpoonup \Psi_\infty$ weakly in $L^2$; 3- we finish
the proof of the theorem through the application of the assumptions
\textbf{A1} through \textbf{A4}.
\subsection{Solutions of the Cauchy problem}\label{ssec:well-posed}
\begin{prop}\label{prop:well1}
Let $\epsilon>0$ and $\Psi_0\in\SSS$. There exists a unique weak
solution $\Psi$ of~\eqref{eq:main}-\eqref{eq:initial} with the the
feedback law $u(t)=u_\epsilon(\Psi(t))$ given by~\eqref{eq:feed},
i.e. $\Psi\in C^0(\RR^+,\SSS)$ $\cap
C^1(\RR^+,H^{-2}(\RR^N,\CC))$,~\eqref{eq:main} holds in
$H^{-2}(\RR^N,\CC)$ for every $t\in\RR^+$ and the
equality~\eqref{eq:initial} holds in $\SSS$.

If, moreover, $\Psi_0\in H^2(\RR^N,\CC)$ and multiplication by
$\mu(x)$ defines a bounded operator over $H^2(\RR^N)$, then $\Psi$
is a strong solution, i.e. $\Psi\in C^0([0,T],H^2$ $(\RR^N,\CC))\cap
C^1([0,T],L^2(\RR^N,\CC))$.
\end{prop}

\textit{Proof.} Let $M\in \NN^*$  be the number of bound states of
$H_0=-\triangle+V(x)$ and $T>0$ such that
\begin{equation}\label{eq:T}
2(M+1)c\|\mu\|_{L^\infty}^2 T e^{(M+1)c\|\mu\|_{L^\infty}^2~T}<1.
\end{equation}
In order to build solution on $[0,T]$, we apply the Banach
fixed-point Theorem to the following map
$$\begin{array}{cccc}
\Theta : & C^{0}([0,T],\SSS)    & \rightarrow & C^{0}([0,T],\SSS) \\
         & \xi                  & \mapsto     & \Psi
\end{array}$$
where $\Psi$ is the solution of~\eqref{eq:main}-\eqref{eq:initial}
with $u(t)=u_\epsilon(\xi(t))$.

The map $\Theta$ is well defined and maps $C^{0}([0,T],\SSS)$ into
itself. Indeed, when $\xi \in C^{0}([0,T],\SSS)$, $u:t \mapsto
u_{\epsilon}(\xi(t))$ is continuous and thus the
Proposition~\ref{prop:well-posed} ensures the existence of a unique
weak solution $\Psi$. Notice that the map $\Theta$ takes values in
$C^{0}([0,T],\SSS)) \cap C^{1}([0,T],H^{-2}_{(0)})$.

Let us prove that $\Theta$ is a contraction of $C^{0}([0,T],\SSS)$.
Let $\xi_{j} \in C^{0}([0,T],\SSS)$, $u_{j}:=u_{\epsilon}(\xi_{j})$,
$\Psi_{j}:=\Theta(\xi_{j})$, for $j=1,2$ and
$\Delta:=\Psi_{1}-\Psi_{2}$. We have
$$
\Delta(t)=i \int_{0}^{t} e^{-i(t-s)H_0} [ u_{1} \mu(x) \Delta (s) +
(u_{1}-u_{2}) \mu(x) \Psi_{2}(s) ]ds.
$$
Thanks to~\eqref{eq:feed}, we have $\|u_{j}\|_{L^{\infty}(0,T)}
\leqslant (M+1)c\|\mu\|_{L^\infty}$ for $j=1,2$ and
$\|v_{1}-v_{2}\|_{L^{\infty}(0,T)} \leqslant
2(M+1)c\|\mu\|_{L^\infty} \|\xi_{1}-\xi_{2}\|_{C^{0}([0,T],L^{2})}$.
Thus
\begin{equation} \label{Delta-intermediaire}
\|\Delta(t)\|_{L^{2}} \leqslant \int_{0}^{t}
(M+1)c\|\mu\|_{L^\infty}^2 \|\Delta(s)\|_{L^{2}} + 2(M+1)c
\|\mu\|_{L^\infty}^2 \|\xi_{1}-\xi_{2}\|_{C^{0}([0,T],L^{2})} ds.
\end{equation}
Therefore, the Gronwall Lemma implies
$$
\|\Delta(t)\|_{C^{0}([0,T],L^{2})} \leqslant 2(M+1)c
\|\mu\|_{L^\infty}^2 T e^{(M+1)c\|\mu\|_{L^\infty}^2
T}~\|\xi_{1}-\xi_{2}\|_{C^{0}([0,T],L^{2})},
$$
and so~\eqref{eq:T} ensures that $\Theta$ is a contraction of the
Banach space $C^{0}([0,T],\SSS)$. Therefore, there exists a fixed
point $\Psi \in C^{0}([0,T],\SSS)$ such that $\Theta(\Psi)=\Psi$.
Since $\Theta$ takes values in $C^{0}([0,T],\SSS) \cap
C^{1}([0,T],H^{-2}(\RR^N,\CC))$, necessarily $\Psi$ belongs to this
space, thus, it is a weak solution
of~\eqref{eq:main}-\eqref{eq:initial} on $[0,T]$.

If, moreover, $\Psi_{0} \in H^{2}(\RR^N,\CC)$ and multiplication by
$\mu(x)$ defines a bounded operator over $H^2(\RR^N,\RR)$, then
applying the Proposition~\ref{prop:well-posed}, the map $\Theta$
takes values in $C^{0}([0,T],H^{2}(\RR^N,\CC)) \cap
C^{1}([0,T],L^2(\RR^N,\CC))$ thus $\Psi$ belongs to this space and
it is a strong solution.

Finally, we have introduced a time $T>0$ and, for every $\Psi_{0}
\in \SSS$, we have built a weak solution $\Psi \in
C^{0}([0,T],\SSS)$ of ~\eqref{eq:main}-\eqref{eq:initial} on
$[0,T]$. Thus, for a given initial condition $\Psi_{0} \in \SSS$, we
can apply this result on $[0,T]$, $[T,2T]$, $[2T,3T]$ etc. This
proves the existence and uniqueness of a global weak solution for
the closed-loop system. $\square$

Note that, by Assumption \textbf{A1}, the initial state $\Psi_0$ is
spanned by the exponentially decaying bound states and therefore
$\Psi_0\in \SSS\cap H^2(\RR^N)$. This, together with the
Proposition~\ref{prop:well1}, terminates the proof of the
well-posedness part of the Theorem~\ref{thm:ggg}.

\subsection{Weak $\omega$-limit set}
Before studying the weak $\omega$-limit set of the closed-loop
system, let us announce two simple and two rather complicated Lemmas
that we will need to characterize this asymptotic regime.

\begin{lem}\label{lem:uL2}
The feedback law $u=u_\epsilon(\Psi)$ defined by~\eqref{eq:feed} is
a member of $L^2_t(\RR^+,\RR)$. In particular, for any $\gamma>0$
there exists $T_\gamma>0$ large enough such that:
$$
\int_{T_\gamma}^\infty |u_\epsilon(\Psi(s))|^2 ds \leq \gamma.
$$
\end{lem}

\textit{Proof.} By definition, we have
$\frac{d\VVV_\epsilon}{dt}=-\frac{1}{c}|u_\epsilon(\Psi)|^2$. The
Lyapunov function $\VVV_\epsilon(\Psi)$ being a decreasing
non-negative function, there exists a positive constant $\nu$ such
that $\VVV_\epsilon(\Psi(t))\searrow\nu \geq 0$. Therefore, we have
$$
\int_0^\infty |u_\epsilon(\Psi(t))|^2 dt=-c\int_0^\infty
\frac{d\VVV_\epsilon}{dt}=c(\VVV_\epsilon(\Psi_0)-\nu)<\infty.
$$
$\square$

\begin{lem}\label{lem:weaklim}
Let $\Psi(t)$ denote the weak (or strong) solution of the
closed-loop system. There exists a sequence of times
$(t_n)_{n=1}^\infty\nearrow\infty$ and some function $\Psi_\infty\in
L^2(\RR^N,\CC)$ (with $\|\Psi_\infty\|_{L^2}\leq 1$) such that:
\begin{alignat}{2}
\Psi(t_n) & \rightharpoonup \Psi_\infty \qquad \text{weakly in }
L^2(\RR^N,\CC),\label{eq:convweak}\\
\PP_{\text{disc}}\Psi(t_n) & \rightarrow
\PP_{\text{disc}}\Psi_\infty \qquad \text{strongly in }
L^2(\RR^N,\CC).\label{eq:convPd}
\end{alignat}
\end{lem}

\textit{Proof.} The solution $\Psi$ belonging to $C^0(\RR^+,\SSS)$,
we have $\|\Psi(t)\|_{L^2}=1$ and therefore the existence of a
subsequence $(t_n)_1^\infty\nearrow\infty$ and $\Psi_\infty\in L^2$
such that~\eqref{eq:convweak} holds true is trivial. Moreover,
$$
\|\Psi_\infty\|_{L^2}\leq \liminf_{n\rightarrow\infty}
\|\Psi(t_n)\|_{L^2}=1.
$$
The key to the proof of~\eqref{eq:convPd}, is in the fact that
$\EEE_{\text{disc}}=\text{Range}(\PP_{\text{disc}})$ is finite
dimensional. Indeed, the weak convergence~\eqref{eq:convweak}
implies
$$
\her{\Psi(t_n)}{\phi_j}\rightarrow \her{\Psi_\infty}{\phi_j},\qquad
j=0,1,\cdots,M,
$$
and therefore
$$
\PP_{\text{disc}}\Psi(t_n)= \sum_{j=0}^M
\her{\Psi(t_n)}{\phi_j}\phi_j\rightarrow \sum_{j=0}^M
\her{\Psi_\infty}{\phi_j}\phi_j=\PP_{\text{disc}}\Psi_\infty\quad
\text{strongly in }L^2.
$$
$\square$

\begin{lem}\label{lem:mulim}
Let $\Psi(t)$ denote the weak (or strong) solution of the
closed-loop system. Consider a sequence of times
$(t_n)_{n=1}^\infty\nearrow\infty$ and some strictly positive time
constant $\tau>0$. We have
\begin{equation}\label{eq:mulim}
S(\tau)\PP_{\text{ac}}\Psi(t_n) \rightarrow 0 \qquad \text{strongly
in } L^2_{\mu^2}(\RR^N,\CC)
\end{equation}
\end{lem}

\textit{Proof.} We have the Duhamel's formula:
$$
S(\tau)\Psi(t_n)=S(t_n+\tau)\Psi_0+\frac{1}{i}\int_0^{t_n}
u_\epsilon(\Psi(s))~S(t_n+\tau-s)\mu(x)\Psi(s)ds.
$$
and therefore
\begin{multline}\label{eq:duhamel1}
\|\mu S(\tau)\PP_{\text{ac}}\Psi(t_n)\|_{L^2_x}\leq \|\mu
S(t_n+\tau)\PP_{\text{ac}}\Psi_0\|_{L^2_x}\\+\int_0^{t_n}
|u_\epsilon(\Psi(s))|~\|\mu(x)~S(t_n+\tau-s)\PP_{\text{ac}}\mu(x)\Psi(s)\|_{L^2_x}ds,
\end{multline}
where we have applied the fact that the semigroup operator $S(t)$ of
$H_0$ commutes with the eigenprojection operator $\PP_{\text{ac}}$
of the same Hamiltonian.

We know by the assumption of the Theorem~\ref{thm:ggg} on $\mu$ that
$\mu\in L^p(\RR^N)$ where $p\in[2,2N)$. Applying the Holder
inequality, we have
\begin{equation}\label{eq:est1}
\|\mu(x)~S(t)\PP_{\text{ac}}\psi\|_{L^2_x}\leq
\|\mu\|_{L^p_x}\|S(t)\PP_{\text{ac}}\psi\|_{L^q_x},\qquad
\frac{1}{p}+\frac{1}{q}=\frac{1}{2},
\end{equation}
where $\psi\in L^2(\RR^N)$.

Moreover, applying the dispersive estimate of the
Corollary~\ref{cor:Lpq}, we have
\begin{equation}\label{eq:est2}
\|S(t)\PP_{\text{ac}}\psi\|_{L^q_x}\leq
|t|^{-\frac{N}{p}}\|\psi\|_{L^{q'}_x},\qquad
\frac{1}{q}+\frac{1}{q'}=1,
\end{equation}
for $\psi\in L^{q'}\cap L^2$.

Let us apply these estimates~\eqref{eq:est1} and~\eqref{eq:est2} to
the inequality~\eqref{eq:duhamel1}. For the first term
in~\eqref{eq:duhamel1}, we have
\begin{equation}\label{eq:aux1}
\|\mu~S(t_n+\tau)\PP_{\text{ac}}\Psi_0\|_{L^2_x}\leq
|t_n+\tau|^{-\frac{N}{p}} \|\mu\|_{L^p_x}\|\Psi_0\|_{L^{q'}_x},
\end{equation}
where we have used the fact that $\Psi_0\in L^{q'}$ as it is a
linear combination of the bound states and therefore decaying
exponentially~\cite{agmon-book}. For the second term, we have
\begin{multline}\label{eq:aux2}
\|\mu~S(t_n+\tau-s)\PP_{\text{ac}}\mu(x)\Psi(s)\|_{L^2_x}\leq
|t_n+\tau-s|^{-\frac{N}{p}} \|\mu\|_{L^p_x}\|\mu\Psi(s)\|_{L^{q'}_x}\\
\leq |t_n+\tau-s|^{-\frac{N}{p}} \|\mu\|^2_{L^p_x}
\|\Psi(s)\|_{L^2_x}= |t_n+\tau-s|^{-\frac{N}{p}} \|\mu\|^2_{L^p_x}.
\end{multline}
Here, to obtain the second line from the first one, we have applied
a holder inequality noting that
$\frac{1}{q'}=\frac{1}{2}+\frac{1}{p}$.

Furthermore, for any $\gamma>0$ taking $t_n>T_\gamma$ (where
$T_\gamma$ is given by Lemma~\ref{lem:uL2}), we have
\begin{multline}\label{eq:aux3}
\int_0^{t_n}
|u_\epsilon(\Psi(s))|~\|\mu(x)~S(t_n+\tau-s)\PP_{\text{ac}}\mu(x)\Psi(s)\|_{L^2_x}ds=\\\int_0^{T_\gamma}
|u_\epsilon(\Psi(s))|~\|\mu(x)~S(t_n+\tau-s)\PP_{\text{ac}}\mu(x)\Psi(s)\|_{L^2_x}ds+\\
 \int_{T_\gamma}^{t_n}
|u_\epsilon(\Psi(s))|~\|\mu(x)~S(t_n+\tau-s)\PP_{\text{ac}}\mu(x)\Psi(s)\|_{L^2_x}ds.
\end{multline}
Inserting the estimate~\eqref{eq:aux2} in the first integral
of~\eqref{eq:aux3}, we have
\begin{multline}\label{eq:aux4}
\int_0^{T_\gamma}
|u_\epsilon(\Psi(s))|~\|\mu(x)~S(t_n+\tau-s)\PP_{\text{ac}}\mu(x)\Psi(s)\|_{L^2_x}ds\leq\\
\|u_\epsilon(\Psi(t))\|_{L^2_t}\|\mu\|^2_{L^p_x}\left(\int_0^{T_\gamma}(t_n+\tau-s)^{-\frac{2N}{p}}ds\right)^{1/2}
\leq\\
\frac{\sqrt
p}{\sqrt{|2N-p|}}\|u_\epsilon(\Psi(t))\|_{L^2_t}\|\mu\|^2_{L^p_x}
\left(|t_n+\tau-T_\gamma|^{-\frac{2N-p}{p}}-|t_n+\tau|^{-\frac{2N-p}{p}}\right)^{1/2},
\end{multline}
where we have applied the Cauchy-Schwartz inequality. Note, in
particular that, $p$ being strictly less than $2N$, $\frac{2N-p}{p}$
is strictly positive, and therefore the above
integral~\eqref{eq:aux4} tends to 0 as $t_n\rightarrow\infty$.

Applying once again the Cauchy-Schwartz inequality, this time for
the second integral in~\eqref{eq:aux3}, we have
\begin{multline}\label{eq:aux5}
\int_{T_\gamma}^{t_n}
|u_\epsilon(\Psi(s))|~\|\mu(x)~S(t_n+\tau-s)\PP_{\text{ac}}\mu(x)\Psi(s)\|_{L^2_x}ds\leq\\
\|\mu\|^2_{L^p_x}
\left(\int_{T_\gamma}^{\infty}|u_\epsilon(\Psi(t))|^2dt\right)^{1/2}
\left(\int_{T_\gamma}^{t_n}(t_n+\tau-s)^{-\frac{2N}{p}}ds\right)^{1/2}\leq\\
\frac{\sqrt p}{\sqrt{|2N-p|}}\gamma^{1/2}\|\mu\|^2_{L^p_x}
\left(|\tau|^{-\frac{2N-p}{p}}-|t_n+\tau-T_\gamma|^{-\frac{2N-p}{p}}\right)^{1/2},
\end{multline}
where, we have used the fact that by definition of $T_\gamma$,
$\int_{T_\gamma}^{\infty}|u_\epsilon(\Psi(t))|^2 dt<\gamma$. In
particular, this implies
\begin{multline}\label{eq:aux6}
\liminf_{t_n\rightarrow \infty} \int_{T_\gamma}^{t_n}
|u_\epsilon(\Psi(s))|~\|\mu(x)~S(t_n+\tau-s)\PP_{\text{ac}}\mu(x)\Psi(s)\|_{L^2_x}ds
\\\leq
\frac{\sqrt
p}{\sqrt{|2N-p|}}\gamma^{1/2}\|\mu\|^2_{L^p_x}|\tau|^{-\frac{2N-p}{2p}}.
\end{multline}
Gathering~\eqref{eq:aux4} and~\eqref{eq:aux6}, we have shown
\begin{multline}\label{eq:aux7}
\liminf_{t_n\rightarrow \infty} \int_{0}^{t_n}
|u_\epsilon(\Psi(s))|~\|\mu(x)~S(t_n+\tau-s)\PP_{\text{ac}}\mu(x)\Psi(s)\|_{L^2_x}ds
\\\leq
\frac{\sqrt
p}{\sqrt{|2N-p|}}\gamma^{1/2}\|\mu\|^2_{L^p_x}|\tau|^{-\frac{2N-p}{2p}}.
\end{multline}
Note, however, that we can choose the constant $\gamma>0$ as small
as we want and therefore we have:
\begin{equation}\label{eq:aux8}
\lim_{t_n\rightarrow \infty} \int_{0}^{t_n}
|u_\epsilon(\Psi(s))|~\|\mu(x)~S(t_n+\tau-s)\PP_{\text{ac}}\mu(x)\Psi(s)\|_{L^2_x}ds=0.
\end{equation}
This, together with~\eqref{eq:aux1}, finishes the proof of
Lemma~\ref{lem:mulim} and we have
\begin{equation}\label{eq:aux9}
\lim_{t_n\rightarrow \infty}\|\mu
S(\tau)\PP_{\text{ac}}\Psi(t_n)\|_{L^2_x}=0.
\end{equation}
$\square$

Applying the above Lemmas, we have the following Lemma, proving the
continuity of the solution of the closed-loop system with respect to
its initial state in the $L^2_{\text{disc}}$-topology.
\begin{lem}\label{lem:cont}
Let $\Psi(t)$ denote the weak (or strong) solution of the
closed-loop system. Consider the time sequence
$\{t_n\}_{n=1}^\infty\nearrow\infty$ and the weak limit
$\Psi_\infty$ as in Lemma~\ref{lem:weaklim} and define
$\Psi_{\infty,\text{disc}}=\PP_{\text{disc}}\Psi_\infty$. Consider
the two closed-loop systems
\begin{alignat}{2}
i\dotex\Psi_n&=-\triangle\Psi_n+V(x)\Psi_n+u_\epsilon(\Psi_n)\mu(x)\Psi_n,\qquad\qquad
\Psi_n|_{t=0}=\Psi(t_n),\label{eq:closed1}\\
i\dotex\widetilde\Psi
&=-\triangle\widetilde\Psi+V(x)\widetilde\Psi+u_\epsilon(\widetilde\Psi)\mu(x)\widetilde\Psi,\qquad\qquad\widetilde\Psi|_{t=0}=\widetilde\Psi_0=\Psi_{\infty,\text{disc}}.\label{eq:closed2}
\end{alignat}
We have, for any $\tau>0$, that
\begin{equation}\label{eq:contPd}
\PP_{\text{disc}}\Psi_n(\tau)\rightarrow
\PP_{\text{disc}}\widetilde\Psi(\tau)\qquad\text{ strongly in } L^2
\text{ as } n\rightarrow \infty.
\end{equation}
\end{lem}

\textit{Proof.} In this aim, we consider a stronger semi-norm than
$L^2_{\text{disc}}$ defined by
$\|\psi\|_{\HHH}=\max(\|\PP_\text{disc}\psi\|_{L^2},\|\mu\psi\|_{L^2})$.
Note however that this semi-norm is weaker than the
$L^2(\RR^N)$-norm
$$
\|\PP_\text{disc}\psi\|_{L^2}\leq
\|\psi\|_{L^2}\qquad\text{and}\qquad \|\mu\psi\|_{L^2}\leq
\|\mu\|_{L^\infty}\|\psi\|_{L^2},
$$
and therefore
\begin{equation}\label{eq:comparenorm}
\|\psi\|_{\HHH}\leq \kappa\|\psi\|_{L^2(\RR^N)},\qquad \text{where }
\kappa=\max(1,\|\mu\|_{L^\infty}).
\end{equation}
It is clear that this is enough to prove
\begin{equation}\label{eq:contH}
\|\Psi_n(\tau)-\widetilde\Psi(\tau)\|_{\HHH}\rightarrow 0\qquad
\text{ as } n\rightarrow \infty.
\end{equation}
We have by the Duhamel's formula
\begin{align*}
\Psi_n(\tau)&=S(\tau)\Psi(t_n)+\frac{1}{i}\int_{0}^\tau
u_\epsilon(\Psi_n(s))S(\tau-s)\mu(x)\Psi_n(s,x)ds\\
\widetilde\Psi(\tau)&=S(\tau)\widetilde\Psi_0+\frac{1}{i}\int_{0}^\tau
u_\epsilon(\widetilde\Psi(s)) S(\tau-s)\mu(x)\widetilde\Psi(s,x)ds,
\end{align*}
Noting by $\delta\Psi_n(\tau)=\Psi_n(\tau)-\widetilde\Psi(\tau)$, we
have
\begin{multline*}
\delta\Psi_n(\tau)= S(\tau)(\Psi(t_n)-\widetilde \Psi_0)+
\frac{1}{i}\int_0^\tau u_\epsilon(\Psi_n(s))S(\tau-s)\mu(x)\delta\Psi_n(s)ds\\
+\frac{1}{i}\int_0^\tau\left[u_\epsilon(\Psi_n(s))-u_\epsilon(\widetilde\Psi(s))\right]S(\tau-s)\mu(x)\widetilde\Psi(s)ds.
\end{multline*}
This implies
\begin{multline}\label{eq:aux10}
\|\delta\Psi_n(\tau)\|_{\HHH}\leq
\|S(\tau)(\Psi(t_n)-\widetilde\Psi_0)\|_{\HHH}\\+ \kappa \int_0^\tau
|u_\epsilon(\Psi_n(s))|\|S(\tau-s)\mu(x)\delta\Psi_n(s)\|_{L^2(\RR^N)}ds\\+
\kappa \int_0^\tau
\left|u_\epsilon(\Psi_n(s))-u_\epsilon(\widetilde\Psi(s))\right|\|S(\tau-s)\mu(x)\widetilde\Psi(s)\|_{L^2(\RR^N)}ds,
\end{multline}
where we have applied the inequality~\eqref{eq:comparenorm}.

Furthermore, noting that $S(t)$ induces an isometry over the space
$L^2(\RR^N,\CC)$, we have
\begin{align*}
\|S(\tau-s)\mu(x)\delta\Psi_n(s)\|_{L^2(\RR^N)}&
=\|\mu(x)\delta\Psi_n(s)\|_{L^2(\RR^N)}\leq
\|\delta\Psi_n(s)\|_{\HHH},\\
\|S(\tau-s)\mu(x)\widetilde\Psi(s)\|_{L^2(\RR^N)}& =
\|\mu(x)\widetilde\Psi(s)\|_{L^2(\RR^N)}\leq
\|\mu\|_{L^\infty}\|\widetilde\Psi(s)\|_{L^2}\leq
\|\mu\|_{L^\infty}.
\end{align*}
where, for the second line, we have also applied
$$
\|\widetilde\Psi(s)\|_{L^2}=\|\widetilde\Psi(0)\|_{L^2}\leq
\|\Psi_\infty\|_{L^2}\leq 1.
$$
Inserting the above inequalities in~\eqref{eq:aux10}, we have
\begin{multline}\label{eq:aux11}
\|\delta\Psi_n(\tau)\|_{\HHH}\leq
\|S(\tau)(\Psi(t_n)-\widetilde\Psi_0)\|_{\HHH}+ \kappa
\|u_\epsilon\|_{L^\infty_t}\int_0^\tau
\|\delta\Psi_n(s)\|_{\HHH}ds\\+ \kappa \|\mu\|_{L^\infty}\int_0^\tau
\left|u_\epsilon(\Psi_n(s))-u_\epsilon(\widetilde\Psi(s))\right|ds.
\end{multline}
Note, in particular that, by the definition of the feedback law
$u_\epsilon$,
$\|u_\epsilon\|_{L^\infty_t}<c(M+1)\|\mu\|_{L^\infty}$. Let us study
the second line of~\eqref{eq:aux11}. We have
\begin{multline*}
\left| u_\epsilon(\Psi_n(s))-u_\epsilon(\widetilde\Psi(s)) \right|\leq\\
c(1-\epsilon)\sum_{j=0}^M \left|\her{\mu\Psi_n(s)}{\phi_j}\her{\phi_j}{\Psi_n(s)}-\her{\mu\widetilde\Psi(s)}{\phi_j}\her{\phi_j}{\widetilde\Psi(s)}\right|\\
+c\epsilon\left|\her{\mu\Psi_n(s)}{\phi_0}\her{\phi_0}{\Psi_n(s)}-\her{\mu\widetilde\Psi(s)}{\phi_0}\her{\phi_0}{\widetilde\Psi(s)}\right|,
\end{multline*}
and for all $j\in\{0,1,...,M\}$
\begin{multline*}
\left|\her{\mu\Psi_n(s)}{\phi_j}\her{\phi_j}{\Psi_n(s)}-\her{\mu\widetilde\Psi(s)}{\phi_j}\her{\phi_j}{\widetilde\Psi(s)}\right|\leq\\
\Big|\her{\mu\delta\Psi_n(s)}{\phi_j}\her{\phi_j}{\Psi_n(s)}\Big|+\left|\her{\mu\widetilde\Psi(s)}{\phi_j}\her{\phi_j}{\delta\Psi_n(s)}\right|\leq\\
\|\mu\delta\Psi_n\|_{L^2(\RR^N)}+\|\mu\|_{L^\infty}
\|\PP_{\text{disc}} \delta\Psi_n(s)\|_{L^2(\RR_N)}\leq
\left(1+\|\mu\|_{L^\infty}\right)\|\delta\Psi_n(s)\|_{\HHH},
\end{multline*}
where, for the last inequality, we have applied the Cauchy-Schwartz
inequality, the facts that $\|\Psi_n(s)\|_{L^2}=1$ and $\|\widetilde
\Psi(s)\|_{L^2}\leq 1$, and that
$|\her{\psi}{\phi_j}|\leq\|\PP_{\text{disc}}\psi\|_{L^2}$.

The above inequality, together with~\eqref{eq:aux11}, implies
\begin{multline}\label{eq:aux12}
\|\delta\Psi_n(\tau)\|_{\HHH}\leq
\|S(\tau)(\Psi(t_n)-\widetilde\Psi_0)\|_{\HHH}\\+ c\kappa
\|\mu\|_{L^\infty}(M+2+\|\mu\|_{L^\infty})\int_0^\tau
\|\delta\Psi_n(s)\|_{\HHH}ds.
\end{multline}
Applying the Gronwall Lemma to~\eqref{eq:aux12}, one only needs to
prove
\begin{equation}\label{eq:contH2}
\|S(\tau)(\Psi(t_n)-\widetilde\Psi_0)\|_{\HHH}\rightarrow
0\qquad\text{ as }n\rightarrow\infty.
\end{equation}
As a first step, we clearly have
\begin{equation}\label{eq:lim1}
\|\PP_{\text{disc}}
S(\tau)(\Psi(t_n)-\widetilde\Psi_0)\|_{L^2}=\|\PP_{\text{disc}}
\Psi(t_n)-\widetilde\Psi_0\|_{L^2}\rightarrow 0,
\end{equation}
where we have used the fact that the semigroup $S(\tau)$ induces an
isometry on $L^2(\RR^N)$, and that the projection operator
$\PP_{\text{disc}}$ commutes with the evolution operator $S(\tau)$.

Moreover applying the fact that,
$\PP_{\text{disc}}+\PP_{\text{ac}}=Id_{L^2}$, we have
$$
\|\mu S(\tau)(\Psi(t_n)-\widetilde \Psi_0)\|_{L^2_x}\leq
\|\mu\|_{L^\infty}\|S(\tau)(\PP_{\text{disc}}\Psi(t_n)-\widetilde\Psi_0)\|_{L^2}+\|\mu
S(\tau)\PP_\text{ac}\Psi(t_n)\|_{L^2}.
$$
Applying~\eqref{eq:lim1}, the first term,
$\|S(\tau)(\PP_{\text{disc}}\Psi(t_n)-\widetilde\Psi_0)\|_{L^2}$
converges toward 0 whenever $n\rightarrow\infty$. Moreover, applying
the Lemma~\ref{lem:mulim}
$$
\|\mu S(\tau)\PP_\text{ac}\Psi(t_n)\|_{L^2}\rightarrow
0\qquad\text{as } n\rightarrow\infty.
$$
Therefore,
\begin{equation}\label{eq:lim2}
\|\mu S(\tau)(\Psi(t_n)-\widetilde\Psi_0)\|_{L^2}\rightarrow
0\qquad\text{as }n\rightarrow\infty.
\end{equation}
The two limits~\eqref{eq:lim1} and~\eqref{eq:lim2} imply the
limit~\eqref{eq:contH2} and therefore finish the proof of the
Lemma~\ref{lem:cont}. $\square$

We are now ready to characterize the weak $\omega$-limit set.
\begin{prop}\label{prop:omega}
Let $\Psi(t)$ denote the weak (or strong) solution of the
closed-loop system. Assume for a sequence
$(t_n)_{n=1}^\infty\nearrow\infty$ of times that
$\Psi(t_n)\rightharpoonup\Psi_\infty\in L^2(\RR^N,\CC)$ weakly in
$L^2(\RR^N,\CC)$ (with $\|\Psi_\infty\|_{L^2}\leq 1$). Define
$\Psi_{\infty,\text{disc}}=\PP_{\text{disc}}\Psi_\infty$. One
necessarily has
\begin{multline*}
u_\epsilon(\Psi_{\infty,\text{disc}})=
c[(1-\epsilon)\sum_{j=0}^M\Im(\her{\mu\Psi_{\infty,\text{disc}}}{\phi_j}\her{\phi_j}{\Psi_{\infty,\text{disc}}})
\\+\epsilon\Im(\her{\mu\Psi_{\infty,\text{disc}}}{\phi_0}\her{\phi_0}{\Psi_{\infty,\text{disc}}})]=0.
\end{multline*}
\end{prop}

\textit{Proof.} Consider the Lyapunov function $\VVV_\epsilon(\Psi)$
defined in~\eqref{eq:lyap}. As it is shown in~\eqref{eq:feedlyap},
the choice~\eqref{eq:feed} of $u_\epsilon(\Psi)$ ensures that the
Lyapunov function $\VVV_\epsilon(\Psi(t))$ is a decreasing function
of time. The Lyapunov function $\VVV_\epsilon$ being a positive
function~\eqref{eq:pos}, we have
\begin{equation}\label{eq:limlyap}
\lim_{t\rightarrow\infty}\VVV_\epsilon(\Psi(t))=\eta,
\end{equation}
where $\eta\geq 0$ is a positive constant.

Consider now, the sequence $\{t_n\}_{n=1}^\infty\nearrow\infty$ of
times. The Lyapunov function $\VVV_\epsilon(\psi)$ is trivially
continuous with respect to $\psi$ for the $L^2$-weak topology.
Therefore, as $\Psi_\infty$ is the weak limit of $\Psi(t_n)$, we
have
$$
\VVV_\epsilon(\Psi_\infty)=\lim_{n\rightarrow\infty}
\VVV_\epsilon(\Psi(t_n))=\eta.
$$
Furthermore, noting that the Lyapunov function $\VVV_\epsilon$ only
deals with the population of the bound states, we have
$$
\VVV_\epsilon(\psi)=\VVV_\epsilon(\PP_{\text{disc}}\psi),
$$
and therefore
\begin{equation}\label{eq:omega1}
\VVV_\epsilon(\Psi_{\infty,\text{disc}})=\eta.
\end{equation}
As in the Lemma~\ref{lem:cont}, let us consider the closed-loop
Schr\"odinger equation with the wavefunction $\widetilde\Psi$ and
the initial state $\widetilde\Psi_0=\Psi_{\infty,\text{disc}}$.
Applying Lemma~\ref{lem:cont}, for any $\tau>0$, we have
$$
\PP_{\text{disc}}\Psi(t_n+\tau)\rightarrow
\widetilde\Psi(\tau)\qquad\text{Strongly in } L^2(\RR^N,\CC) \text{
as }n\rightarrow\infty.
$$
As the Lyapunov function $\VVV_\epsilon(\psi)$ is continuous with
respect to $\psi$ for the $L^2_{\text{disc}}$ semi-norm, we have
$$
\VVV_\epsilon(\Psi(t_n+\tau))\rightarrow\VVV_\epsilon(\widetilde\Psi(\tau))\qquad\text{as
}n\rightarrow\infty.
$$
But, applying~\eqref{eq:limlyap}, we know that
$$
\VVV_\epsilon(\Psi(t_n+\tau))\rightarrow\eta\qquad\text{as
}n\rightarrow\infty,
$$
and therefore,
\begin{equation}\label{eq:omega2}
\VVV_\epsilon(\widetilde\Psi(\tau))=\eta=\VVV_\epsilon(\Psi_{\infty,\text{disc}})=\VVV_\epsilon(\widetilde\Psi(0)).
\end{equation}
Thus, the Lyapunov function $\VVV_\epsilon$ remains constant on the
closed-loop trajectory of $\widetilde\Psi(t)$. This, together
with~\eqref{eq:feedlyap} and~\eqref{eq:feed}, implies
\begin{equation}\label{eq:omega3}
\frac{\partial}{\partial
\tau}\VVV_\epsilon(\widetilde\Psi(\tau))=-\frac{1}{c}u_\epsilon^2(\widetilde\Psi(\tau))=0,
\end{equation}
and therefore by continuity of $u_\epsilon(\widetilde\Psi(\tau))$
with respect to $\tau$ and passing to the limit at $\tau=0$, we can
finish the proof of the Proposition~\ref{prop:omega}. $\square$

\subsection{Non-degeneracy assumptions and the proof of
Theorem~\ref{thm:ggg}}\label{ssec:thmggg}

We have now all the elements to finish the proof of the
Theorem~\ref{thm:ggg}.
\begin{prop}\label{prop:A3A4}
Let $\Psi(t)$ denote the weak (or strong) solution of the
closed-loop system. Consider the sequence
$(t_n)_{n=1}^\infty\nearrow\infty$, the weak limit $\Psi_\infty$,
and its discrete part $\Psi_{\infty,\text{disc}}$ as in
Proposition~\ref{prop:omega}. Under the assumptions \textbf{A1}
through \textbf{A4} of Theorem~\ref{thm:ggg}, we have
\begin{equation}\label{eq:limPd}
\Psi_{\infty,\text{disc}}=\varsigma \phi_0,\qquad
|\varsigma|^2>1-\epsilon.
\end{equation}
\end{prop}

\textit{Proof.} Define $\eta$ as in~\eqref{eq:limlyap}. We now, in
particular that,
\begin{equation}\label{eq:varrho}
\VVV_\epsilon(\Psi_\infty)=\eta\leq \VVV_\epsilon(\Psi(0)) <\epsilon
\end{equation}
where we have applied~\eqref{eq:initlyap} (and therefore the
assumptions \textbf{A1} and \textbf{A2}).

Let us take
$$
\Psi_{\infty,\text{disc}}=\sum_{j=0}^M \varsigma_j \phi_j.
$$
Taking the closed-loop system $\widetilde\Psi(t)$ as in the proof of
the Proposition~\ref{prop:omega}, we have by~\eqref{eq:omega3} that
$u_\epsilon(\widetilde{\Psi}(\tau))=0$. Therefore the wavefunction
$\widetilde\Psi(\tau)$ evolves freely with the Hamiltonian
$H_0=-\triangle+V(x)$ and so is given as follows
$$
\widetilde\Psi(\tau)=\sum_{j=0}^M\varsigma_j e^{-i\lambda_j
\tau}\phi_j.
$$
By~\eqref{eq:omega3} we have
\begin{multline*}
u(\widetilde\Psi)=c(1-\epsilon)\sum_{j,k=0}^M\bar\varsigma_j\varsigma_k
e^{i(\lambda_j-\lambda_k)\tau}\her{\mu\phi_k}{\phi_j}
\\+c\epsilon\sum_{j=0}^M\bar\varsigma_0\varsigma_j
e^{i(\lambda_0-\lambda_j)\tau}\her{\mu\phi_j}{\phi_0}=0\qquad\forall
t\geq 0.
\end{multline*}
The assumption \textbf{A3} of non-degenerate transitions applies
now. As the above relation holds true for any $\tau\geq 0$, we can
easily see that
$$
\bar\varsigma_j\varsigma_k\her{\mu\phi_k}{\phi_j}=0\qquad \forall
j,k\in\{0,1,...,M\}.
$$
This together with the assumption~\textbf{A4} of simple couplings
imply
$$
\bar\varsigma_j\varsigma_k=0\qquad\forall j,k\in\{0,1,...,M\}.
$$
Thus
$$
\exists j\in\{0,1,...,M\} \quad\text{such that }
\varsigma_j=\varsigma\neq0 \quad\text{and }\varsigma_k=0\quad
\forall k\neq j.
$$
We show that the only possibility for this index $j$ is to be 0. If
this is not the case ($j\neq 0$) taking
$\widetilde\Psi=\varsigma\phi_j$ with $|\varsigma|\leq 1$,
$$
\eta=\VVV_\epsilon(\Psi_\infty)=1-(1-\epsilon)|\varsigma|^2\geq
\epsilon,
$$
which is obviously in contradiction with~\eqref{eq:varrho}. Thus
$$
\Psi_{\infty,\text{disc}}=\varsigma\phi_0
$$
with $|\varsigma|\leq 1$. Therefore
$$
\VVV_\epsilon(\Psi_\infty)=1-(1-\epsilon)|\varsigma|^2-\epsilon|\varsigma|^2
=1-|\varsigma|^2.
$$
Apply once again~\eqref{eq:varrho}, we have
$1-|\varsigma|^2<\epsilon$ and so we can finish the proof of the
Proposition~\ref{prop:A3A4}. $\square$

Let us now, finish the proof of Theorem~\ref{thm:ggg}.

\textit{Proof of Theorem~\ref{thm:ggg}} The well-posedness of the
closed-loop system has been addressed in
Proposition~\ref{prop:well1}. In order to prove the approximate
stabilization result, let us assume that there exists a sequence of
times $\{\tilde t_n\}_{n=0}^\infty\nearrow\infty$ such that
\begin{equation}\label{eq:absurd}
|\her{\Psi(\tilde t_n)}{\phi_0}|^2\leq 1-\epsilon \qquad\forall n.
\end{equation}
As in Proposition~\ref{prop:A3A4}, we can extract from this sequence
a subsequence (noted still by $\{\tilde t_n\}_{n=0}^\infty$ for
simplicity sakes) such that
$$
\PP_{\text{disc}}\Psi(\tilde
t_n)\stackrel{L^2-\text{strong}}\longrightarrow \widetilde
\varsigma\phi_0\quad\text{as } n\rightarrow\infty,
$$
with $|\widetilde\varsigma|^2>1-\epsilon$. This obviously implies
$$
\liminf_{t\rightarrow\infty}|\her{\Psi(\tilde
t_n)}{\phi_0}|^2>1-\epsilon
$$
and is in contradiction with~\eqref{eq:absurd}. We have therefore
finished the proof of the Theorem~\ref{thm:ggg}. $\square$
\section{Proof of the Theorem~\ref{thm:gas}}\label{sec:proof2}
Let us now get back to the Theorem~\ref{thm:gas}. Comparing to the
Theorem~\ref{thm:ggg}, the only difference is in the fact that, we
are also considering the 1D case and that the interaction
Hamiltonian $\mu\in\LLL\cap L^\infty$ instead of $L^{2N-}\cap
L^\infty$ for the Theorem~\ref{thm:ggg}. Therefore the only cases
remaining to be treated are either the 1D case or the cases where
$\mu\in L^p(\RR^N)\cap L^\infty(\RR^N)$ with $p\geq 2N$.

Considering these cases and following the same steps as in the proof
of the Theorem~\ref{thm:ggg}, the only place where we will have a
problem to proceed the proof is the passage from~\eqref{eq:aux4}
and~\eqref{eq:aux5} to~\eqref{eq:aux6} and~\eqref{eq:aux7}. Indeed,
as $2N-p$ is not strictly positive, we can not ensure the
convergence towards 0 of the terms in~\eqref{eq:aux4}
and~\eqref{eq:aux5}.

A deep study of the estimates~\eqref{eq:aux4} and~\eqref{eq:aux5}
shows that they can be improved if one had $u_\epsilon\in
L^{1+\delta}_t$ for $\delta\in(0,1)$ instead of $L^2_t$ as in the
proof of the Theorem~\ref{thm:ggg}. Indeed, if one could show that
\begin{equation}\label{eq:u}
u_\epsilon\in L^{(\frac{p}{p-N}-\varpi)}_t\qquad\text{for some }
\varpi>0,
\end{equation}
we could replace the estimates~\eqref{eq:aux4} and~\eqref{eq:aux5}
with
\begin{multline}\label{eq:aux4bis}
\int_0^{T_\gamma}
|u_\epsilon(\Psi(s))|~\|\mu(x)~S(t_n+\tau-s)\PP_{\text{ac}}\mu(x)\Psi(s)\|_{L^2_x}ds\leq\\
\|u_\epsilon(\Psi(t))\|_{L^{(\frac{p}{p-N}-\varpi)}_t}\|\mu\|^2_{L^p_x}\left(\int_0^{T_\gamma}(t_n+\tau-s)^{-\frac{\zeta
N}{p}}ds\right)^{1/2}
\leq\\
\frac{\sqrt p}{\sqrt{|\zeta
N-p|}}\|u_\epsilon(\Psi(t))\|_{L^{(\frac{p}{p-N}-\varpi)}_t}\|\mu\|^2_{L^p_x}
\left(|t_n+\tau-T_\gamma|^{-\frac{\zeta
N-p}{p}}-|t_n+\tau|^{-\frac{\zeta N-p}{p}}\right)^{1/2},
\end{multline}
and
\begin{multline}\label{eq:aux5bis}
\int_{T_\gamma}^{t_n}
|u_\epsilon(\Psi(s))|~\|\mu(x)~S(t_n+\tau-s)\PP_{\text{ac}}\mu(x)\Psi(s)\|_{L^2_x}ds\leq\\
\|\mu\|^2_{L^p_x}
\left(\int_{T_\gamma}^{\infty}|u_\epsilon(\Psi(t))|^{\frac{p}{p-N}-\varpi}dt\right)^{1/2}
\left(\int_{T_\gamma}^{t_n}(t_n+\tau-s)^{-\frac{\zeta N}{p}}ds\right)^{1/2}\leq\\
\frac{\sqrt p}{\sqrt{|\zeta N-p|}}\gamma^{1/2}\|\mu\|^2_{L^p_x}
\left(|\tau|^{-\frac{\zeta
N-p}{p}}-|t_n+\tau-T_\gamma|^{-\frac{\zeta N-p}{p}}\right)^{1/2},
\end{multline}
where
\begin{equation}\label{eq:zeta}
\zeta=\frac{p-\varpi(p-N)}{N-\varpi(p-N)},
\end{equation}
noting that we have applied the Holder inequality and the fact that
$$
\frac{1}{\zeta}+\frac{1}{\frac{p}{p-N}-\varpi}=1.
$$
Note that, as $p\geq 2N$ (this is also true for the 1D case as
$p\geq 2$), $\frac{p}{p-N}\in (1,2]$ and therefore there exists some
positive $\varpi>0$ such that $ \frac{p}{p-N}-\varpi>1$.
Furthermore, we have
$$
\zeta=\frac{p-\varpi(p-N)}{N-\varpi(p-N)}>\frac{p}{N},
$$
and therefore $\zeta N-p$ is positive. We can thus proceed the proof
of the Theorem~\ref{thm:gas} following the same steps as those of
the Theorem~\ref{thm:ggg}.

However, it seems that one can not hope to prove an estimate of the
form~\eqref{eq:u} for the feedback law $u_\epsilon$
of~\eqref{eq:feed}. We, therefore, need to change the feedback
strategy. This might be done applying the feedback
law~\eqref{eq:feed2}.

The following Proposition clearly implies the Theorem~\ref{thm:gas}.
\begin{prop}\label{prop:feed2}
Consider the Schr{\"o}dinger
equation~\eqref{eq:main}-~\eqref{eq:initial}. We suppose the
assumptions of the Theorem~\ref{thm:gas} on the potential $V(x)$ and
we take $\mu\in L^p(\RR^N)\cap L^\infty(\RR^N)$ for some $p\geq 2N$.
We suppose moreover the assumptions \textbf{A1} through \textbf{A4}
to hold true.

Then for any $\epsilon >0$, applying the feedback law
$u(t)=u_{\epsilon,\alpha}(\Psi(t))$ of~\eqref{eq:feed2} with
$$
\alpha=\frac{p-2N+\varpi(p-N)}{N-\varpi(p-N)},\qquad
0<\varpi<\frac{N}{N-p},
$$
the closed-loop system admits a unique weak solution in
$C^0([0,T],\SSS)\cap C^1([0,T],$ $H^{-2}(\RR^N,\CC))$. Moreover the
state of the system ends up reaching a population more than
$(1-\epsilon)$ in the eigenfunction $\phi_0$ (approximate
stabilization):
$$
\liminf_{t\rightarrow
\infty}|\her{\Psi(t,x)}{\phi_0(x)}|^2>1-\epsilon.
$$
If, moreover multiplication by $\mu(x)$ defines a bounded operator
over $H^2(\RR^N)$, then $\Psi$ is a strong solution, i.e. $\Psi\in
C^0([0,T],H^2$ $(\RR^N,\CC))\cap C^1([0,T],L^2(\RR^N,\CC))$.
\end{prop}

\textit{Proof.} Considering the Lyapunov function $\VVV_\epsilon$
of~\eqref{eq:lyap}, the choice of the feedback law implies
$$
\frac{d}{dt}\VVV_\epsilon=-\frac{c^{1+\alpha}}{c^{2+\alpha}}|u_{\epsilon,\alpha}|^{\frac{2+\alpha}{1+\alpha}},
$$
and therefore proceeding as in the proof of the Lemma~\ref{lem:uL2}
\begin{equation}\label{eq:uu}
u_{\epsilon,\alpha}\in
L_t^{\frac{2+\alpha}{1+\alpha}}=L^{\frac{p}{p-N}-\varpi}_t.
\end{equation}
In particular, for $\gamma>0$, we will chose $T_\gamma$ such that
\begin{equation}\label{eq:tg}
\int_{T_\gamma}^\infty
|u_\epsilon(\Psi(s))|^{\frac{p}{p-N}-\varpi}ds \leq \gamma.
\end{equation}
One can then proceed the proof of the Proposition~\ref{prop:A3A4},
exactly as in the proof of the Theorem~\ref{thm:ggg}, replacing only
the Lemma~\ref{lem:uL2} with~\eqref{eq:uu} and~\eqref{eq:tg} and the
estimates~\eqref{eq:aux4} and~\eqref{eq:aux5} by~\eqref{eq:aux4bis}
and~\eqref{eq:aux5bis}. $\square$

\section{Relaxations}\label{sec:relax}

As it has been proved in previous sections, the approximate
stabilization of a quantum particle around the bound states of a
decaying potential (satisfying the decay assumption \textbf{(A)})
may be investigated through explicit feedback laws~\eqref{eq:feed}
or~\eqref{eq:feed2}. The assumptions on the potential $V$ or the
interaction Hamiltonian $\mu$ are not so restrictive and seem to be
satisfied for a large class of physical systems. However, the
assumptions~\textbf{A1} through \textbf{A4} may seem to be too
restrictive. In particular, the assumption \textbf{A1} does not
allow the approximate stabilization of an initial  wavefunction with
a non-zero population in the absolutely continuous
part~$\EEE_{\text{ac}}$.

The aim of this Section is to give some ideas to relax these
assumptions and to consider some more general situations. Some
discussions on the assumption \textbf{A1} will be addressed in
subsection~\ref{ssec:A1}. Furthermore, a significant relaxation of
the assumptions \textbf{A3} and \textbf{A4} will be addressed in
subsection~\ref{ssec:A3A4}. Concerning the assumption \textbf{A2},
we only give the following remark which states that this assumption
is, actually, not at all restrictive in practice.
\begin{rem}\label{rem:A2}
Physically, the assumption \textbf{A2} in not really restrictive.
Indeed, even if $\her{\Psi_0}{\phi_0}=0$, a control field in
resonance with the natural frequencies of the system (the difference
between the eigenvalues corresponding to an eigenfunction whose
population in the initial state is non-zero and the ground state)
will, instantaneously, ensure a non-zero population of the ground
state in the wavefunction. Then, one can just apply the feedback law
of the Theorem~\ref{thm:ggg} or~\ref{thm:gas}.
\end{rem}

\subsection{Assumption \textbf{A1}}\label{ssec:A1}

Before discussing an idea which may result in a significant
relaxation of this assumption, let us provide a remark which states
that the result of the Theorem~\ref{thm:gas} still holds true if we
relax slightly the assumption~\textbf{A1}.

\begin{rem}\label{rem:A1}
Consider the Schr{\"o}dinger
equation~\eqref{eq:main}-~\eqref{eq:initial} with the same
assumptions on $V$ and $\mu$ as in Theorem~\ref{thm:gas}. We
consider moreover the assumptions \textbf{A2} through \textbf{A4}
and we replace \textbf{A1} with
\begin{description}
  \item[\textbf{A1'}] $\|\PP_{\text{ac}}\Psi_0\|_{L^2(\RR^N)}<\frac{\epsilon}{1-\epsilon}|\her{\Psi}{\phi_0}|^2$.
\end{description}
The feedback law~\eqref{eq:feed2} still ensures the approximate
stabilization of the closed-loop.

one only needs to note that,
\begin{equation*}
\VVV_{\epsilon}(\Psi_0)<
1-(1-\epsilon)(1-\frac{\epsilon}{1-\epsilon}|\her{\Psi_0}{\phi_0}|^2)-\epsilon|\her{\Psi_0}{\phi_0}|^2=\epsilon,
\end{equation*}
where we have applied the assumptions \textbf{A1'} and \textbf{A2}.
The rest of the proof follows exactly as in the
Theorem~\ref{thm:ggg}.

Here, we have relaxed the assumption~\textbf{A1} by allowing a very
small part of the population of the initial state to belong to the
continuum. In fact, this allowed continuum population is bounded by
an $O(\epsilon)$-proportion of the population in the target state
$\phi_0$.
\end{rem}

The assumption \textbf{A1'} of the Remark~\ref{rem:A1} is still
quite restrictive. The question is therefore to provide a strategy
permitting us to approximately stabilize an important part of the
continuum. Note that the controllability of this particular problem
has never been treated. It seems that one can not in general hope to
have a strong controllability result. In fact, considering the
potentials $V$ and  $\mu$ of compact supports and taking an initial
state of support outside $\text{supp}(V)\cup \text{supp}(\mu)$, it
seems that an important part of the population may be lost at
infinity through the dispersion phenomena and this before the
controller even has the time to see and to influence the state.
However, one might be interested to control a part of the continuum.

Consider for example the potentials $V$ and $\mu$ to be negative and
of compact supports and moreover that $\text{supp}(\mu)\subset
\text{supp}(V)$. Considering the Hamiltonian
$H_\lambda=-\triangle+V+\lambda \mu$ in the strong coupling limit
($\lambda\rightarrow\infty$), this Hamiltonian admits more and more
bound states. One can therefore cover a higher and higher
dimensional subspace of $L^2(\RR^N)$ through the discrete eigenspace
of $H_\lambda$. Assume an initial state $\Psi_0$ which has a large
population in the continuum of $H_0=-\triangle+V$ but a small
population in the continuum of $H_\lambda$ for some $\lambda>0$.
Applying the strategy of the Theorem~\ref{thm:gas}, to the free
Hamiltonian $H_\lambda$ and the interaction Hamiltonian $\mu$ one
may hope to reach an $\epsilon$-neighborhood of an arbitrary bound
state of $H_\lambda$. Note, in particular that while reaching this
bound state the control field $u(t)$ has converged towards
$-\lambda$. Letting now the control field $u(t)\sim-\lambda$ varying
slowly towards zero and applying the quantum adiabatic theory (see
e.g. ~\cite{avron-99}) the state of the system will follow closely a
bound state of the Hamiltonian $H_\lambda$ ($\lambda\rightarrow 0$).
If the target bound state of $H_\lambda$ is chosen to be on the
analytic branch corresponding to the evolution of the ground state
of $H_0$ (see e.g.~\cite{kato-book-66}), as the control tends to 0,
we may reach the $\epsilon$-neighborhood of the desired target
$\phi_0$. This idea of applying the quantum adiabatic theory and the
large coupling limit to ensure the control of a population in the
continuum will be explored in future works.

\subsection{Assumptions \textbf{A3} and
\textbf{A4}}\label{ssec:A3A4}

After the above discussions on the assumptions \textbf{A1} and
\textbf{A2}, let us study the non-degeneracy assumptions \textbf{A3}
and \textbf{A4}. Similar assumptions to \textbf{A3} and \textbf{A4}
have already been considered for the stabilization of finite
dimensional quantum
systems~\cite{mirrahimi-et-al2-04,beauchard-et-al07}. When dealing
with finite dimensional systems, the assumptions~\textbf{A3} and
\textbf{A4} are equivalent to the controllability of the linearized
system around various eigenstates of the system. For these finite
dimensional systems, it was shown, in~\cite{mirrahimi-et-al2-04}
through the quantum adiabatic theory, and
in~\cite{beauchard-et-al07} through the implicit Lyapunov control
techniques, that the non-degeneracy assumptions can be relaxed
significantly. These relaxations have even been applied to the
problem of the approximate stabilization of the quantum particle in
an infinite potential well, being an infinite dimensional system
(see~\cite{beauchard-mirrahimi-07}). In this subsection, we will see
that such relaxations may also be considered for our control problem
of the quantum particle in a decaying potential.

In this aim, we consider the potential $V$ and the interaction
Hamiltonian $\mu$ both to satisfy the decay assumption \textbf{(A)}.
Similarly to the previous works, we consider the family of the
perturbed Hamiltonian $H_\sigma=-\triangle+V(x)+\sigma\mu(x)$ with
$|\sigma|\ll 1$ a small real constant. The family $\{H_\sigma\}$ is
a self-adjoint holomorphic family of type (A) in the sense of Kato
(see~\cite{kato-book-66}, page 375). Thus, the eigenvalues and the
bound states of $H_\sigma$ are holomorphic functions of $\sigma$
around zero.

The absence of zero energy eigenstate for $H_0=-\triangle+V$ as been
assumed in the decay assumption \textbf{(A)} implies the existence
of a strictly positive threshold $\sigma^*$ such that for $\sigma$
evolving in $(-\sigma^*,\sigma^*)$ the bound states $\phi_j$ of
$H_0$ stay bound states of $H_\sigma$ and do not join the continuum.
This ensures that that, the perturbed eigenvalues
$\{\lambda_{\sigma,j}\}_{j=0}^M$ of $H_\sigma$ (with
$\lambda_{0,j}=\lambda_j$) are well-defined and remain less than
zero. Note that, one might have the appearance of new bound states
but the bound states of $H_0$ will not disappear while considering
perturbations of amplitude $|\sigma|<\sigma^*$. We have therefore
the following Theorem:
\begin{thm}\label{thm:A3A4}
Consider the Schr{\"o}dinger
equation~\eqref{eq:main}-~\eqref{eq:initial} with the decay
assumption \textbf{(A)} on both $V$ and $\mu$. We assume moreover
that the space dimension $N\geq 2$ and $\mu\in L^{2N-}\cap
L^\infty$. We consider the assumptions \textbf{A1} and \textbf{A2}
and we replace \textbf{A3} and \textbf{A4} with:
\begin{description}
  \item[\textbf{(A3-A4)'}] there exists $\bar \sigma\in(0,\sigma^*)$
  such that the non-degeneracy assumptions \textbf{A3} and
  \textbf{A4} hold for the eigenvalues and the eigenstates of the
  perturbed Hamiltonian $H_{\bar\sigma}=-\triangle+V+\bar\sigma \mu$.
\end{description}
There exists then a feedback law $u(\Psi)$, such that the
closed-loop system admits a unique weak solution and that
$$
\liminf_{t\rightarrow\infty} |\her{\Psi(t)}{\phi_0}|^2>1-\epsilon.
$$
\end{thm}
\begin{rem}\label{rem:A3A4}
Roughly speaking the non-degeneracy assumptions \textbf{A3} and
\textbf{A4} are always satisfied unless some kind of symmetry is
admitted in the potential. Formally, the assumption
\textbf{(A3-A4)'} states that if we can break this symmetry through
the addition of the interaction Hamiltonian, we are still able to
ensure the approximate stabilization result.
\end{rem}
\begin{rem}
Applying the same technique as in Proposition~\ref{prop:feed2}, one
can extend the result of the Theorem~\ref{thm:A3A4} to the case of
dimensions $N\geq 1$ and $\mu\in\LLL\cap L^\infty$.
\end{rem}

\textit{Proof.} By the analyticity of the eigenvalues
$\lambda_{\sigma,j}$ and $\phi_{\sigma,j}$ with respect to $\sigma$
around zero, the assumption \textbf{(A3-A4)'} ensures the existence
of a strictly positive constant $\sigma^\sharp\in(0,\sigma^*)$ such
that the non-degeneracy assumptions \textbf{A3} and \textbf{A4} hold
true for the perturbed Hamiltonians $H_\sigma$ with any $\sigma$ in
the interval $(0,\sigma^\sharp)$.

Applying once more the analyticity of the bound states
$\{\phi_{j,\sigma}\}_{j=0}^M$ with respect to $\sigma$ in
$(0,\sigma^\sharp)$ implies that one can choose
$\sigma^{\sharp\sharp}\in(0,\sigma^\sharp)$ such that
\begin{multline}\label{eq:sigma}
\|\phi_{j,\sigma}-\phi_j\|_{L^2}<\min(\frac{\epsilon}{4},\frac{\epsilon|\her{\Psi_0}{\phi_0}|^2}{2(M+1)(2-\epsilon)+2\epsilon})\\
\forall j=0,1,\cdots, M\text{ and }  \forall
\sigma\in(0,\sigma^{\sharp\sharp}).
\end{multline}
Now, applying the inequality~\eqref{eq:sigma}, the assumption
\textbf{A2} implies:
\begin{multline}\label{eq:insigma}
|\her{\Psi_0}{\phi_{0,\sigma}}|^2\geq
|\her{\Psi_0}{\phi_{0}}|^2-2\left|\her{\Psi_0}{\phi_{0,\sigma}-\phi_0}\right|\\
\geq |\her{\Psi_0}{\phi_{0}}|^2
-2\|\phi_{j,\sigma}-\phi_j\|_{L^2}>0,\quad\forall
\sigma\in(0,\sigma^{\sharp\sharp}),
\end{multline}
where we have applied the Cauchy-Schwartz inequality and the fact
that
$$
\frac{\epsilon}{2(M+1)(2-\epsilon)+2\epsilon}<\frac{1}{2}.
$$
Let us consider the Schr\"odinger equation $(\Sigma_\sigma)$
characterized by the free Hamiltonian
$H_\sigma=-\triangle+V+\sigma\mu$ and the interaction Hamiltonian
$\mu$ for some $\sigma\in(0,\sigma^{\sharp\sharp})$. Applying the
assumption \textbf{(A3-A4)'}, the assumptions \textbf{A3} and
\textbf{A4} hold true for this system. Moreover, the
inequality~\eqref{eq:insigma} implies the assumption \textbf{A2} for
this system. Finally, applying~\eqref{eq:sigma}, we have
($M_\sigma+1\geq M+1$ is the number of the bound states of
$H_\sigma$)
\begin{align}\label{eq:sigA1}
\PP_{\text{ac}}(H_\sigma)(\Psi_0)&=1-\sum_{j=0}^{M_\sigma}|\her{\Psi_0}{\phi_{j,\sigma}}|^2\leq
1-\sum_{j=0}^{M}|\her{\Psi_0}{\phi_{j,\sigma}}|^2\notag\\
&\leq 1-\sum_{j=0}^M
\left(|\her{\Psi_0}{\phi_{j,\sigma}}|^2-2\|\phi_{j,\sigma}-\phi_j\|_{L^2}\right)\notag\\
&<\frac{2(M+1)\epsilon|\her{\Psi_0}{\phi_0}|^2}{2(M+1)(2-\epsilon)+2\epsilon}<\frac{\epsilon/2}{1-\epsilon/2}|\her{\Psi_0}{\phi_0}|^2.
\end{align}
This implies the assumption \textbf{A1'} or the system
$(\Sigma_\sigma$)(having replaced $\epsilon$ by $\epsilon/2$).

Considering therefore the feedback law
\begin{multline}\label{eq:feed3}
v_{\frac{\epsilon}{2},\sigma}(\Psi(t)):=c(1-\frac{\epsilon}{2})\sum_{j=0}^{M_\sigma}
\her{\mu\Psi(t)}{\phi_{j,\sigma}}\her{\phi_{j,\sigma}}{\Psi(t)}\\
+c\frac{\epsilon}{2}\her{\mu\Psi(t)}{\phi_{0,\sigma}}\her{\phi_{0,\sigma}}{\Psi(t)},\qquad
c>0,
\end{multline}
and applying the Theorem~\ref{thm:ggg}, we ensure the approximate
stabilization result:
\begin{equation}\label{eq:limsig}
\liminf_{t\rightarrow\infty} |\her{\Psi(t)}{\phi_{0,\sigma}}|^2\geq
1-\epsilon/2.
\end{equation}
Note that, the feedback law~\eqref{eq:feed3} means the application
of the feedback $u:=-\sigma+v_{\frac{\epsilon}{2},\sigma}(\Psi(t))$
for the main Schr\"odinger
equation~\eqref{eq:main}-\eqref{eq:initial}. Finally the
limit~\eqref{eq:limsig} implies
\begin{multline}
\liminf_{t\rightarrow\infty} |\her{\Psi(t)}{\phi_{0}}|^2\geq
\liminf_{t\rightarrow\infty}
\left(|\her{\Psi(t)}{\phi_{0,\sigma}}|^2-2\|\phi_{0,\sigma}-\phi_0\|_{L^2}\right)\\\geq
1-\epsilon/2-\epsilon/2=1-\epsilon,
\end{multline}
where once again we have applied~\eqref{eq:sigma}. $\square$

\section{Appendix}

This appendix is devoted to the proof of the
Proposition~\ref{prop:well-posed}.

\textit{Proof of Proposition~\ref{prop:well-posed} } Let $\Psi_{0}
\in \SSS$, $T>0$ and $u \in C^{0}([0,T],\RR)$. Let $T_1 \in (0,T)$
be such that
\begin{equation} \label{u-contract}
\|\mu\|_{L^\infty}\|u\|_{L^{1}(0,T_1)}<1.
\end{equation}
We prove the existence of $\Psi \in C^{0}([0,T_1],L^{2}(\RR^N,\CC))$
such that~\eqref{eq:mild} holds by applying the Banach fixed point
theorem to the map
$$\begin{array}{cccc}
\Theta : & C^{0}([0,T_1],L^{2}) & \rightarrow & C^{0}([0,T_1],L^{2})
\\
            &   \xi              & \mapsto     & \Psi
\end{array}$$
where $\Psi$ is the weak solution of
$$
i \frac{\partial \Psi}{\partial t} =  H_0\Psi - u(t)\mu(x) \xi ,
\qquad \Psi(0,x)=\Psi_{0}(x)
$$
i.e. $\Psi \in C^{0}([0,T_1],L^{2})$ and satisfies, for every $t \in
[0,T_1]$,
$$
\Psi(t)=e^{-iH_0t} \Psi_{0} +i\int_{0}^{t} e^{-iH_0(t-s)} u(s)
\mu(x) \xi(s) ds \text{  in  } L^{2}(\RR^N,\CC).
$$
Notice that $\Theta$ takes values in
$C^{1}([0,T_1],H^{-2}(\RR^N,\CC))$.

For $\xi_{1}, \xi_{2} \in C^{0}([0,T_1],L^{2}(\RR^N,\CC))$,
$\Psi_{1}:=\Theta(\xi_{1})$, $\Psi_{2}:=\Theta(\xi_{2})$ we have
$$
(\Psi_{1}-\Psi_{2})(t)=i\int_{0}^{t} e^{-iH_0(t-s)} u(s)
\mu(x)(\xi_{1}-\xi_{2})(s) ds
$$
thus
$$
\|(\Psi_{1}-\Psi_{2})(t)\|_{L^{2}} \leqslant \|\mu\|_{L^\infty}
\int_{0}^{t} |u(s)| ds \|\xi_{1}-\xi_{2}\|_{C^{0}([0,T_1],L^{2})}.
$$
The assumption~\eqref{u-contract} guarantees that $\Theta$ is a
contraction of $C^{0}([0,T_1],L^{2})$, thus, $\Theta$ has a fixed
point $\Psi \in C^{0}([0,T_1],L^{2})$. Since $\Theta$ takes values
in $C^{1}([0,T_1],H^{-2})$, then $\Psi$ belongs to this space.
Moreover, this function satisfies~\eqref{eq:mild}.

Finally, we have built weak solutions on $[0,T_1]$ for every
$\Psi_{0}$, and the time $T_1$ does not depend on $\Psi_{0}$, thus,
this gives solutions on $[0,T]$.

Let us prove that this solution is continuous with respect to the
the initial condition $\Psi_{0}$, for the
$L^{2}(\RR^N,\CC)$-topology. Let $\Psi_{0}, \Phi_{0} \in \SSS$ and
$\Psi$, $\Phi$ the associated weak solutions. We have
$$
\|(\Psi-\Phi)(t)\|_{L^{2}} \leqslant \| \Psi_{0} - \Phi_{0}
\|_{L^{2}} + \|\mu\|_{L^\infty}\int_{0}^{t} |u(s)|
\|(\Psi-\Phi)(s)\|_{L^{2}} ds,
$$
thus Gronwall Lemma gives
$$
\|(\Psi-\Phi)(t)\|_{L^{2}} \leqslant \| \Psi_{0} - \Phi_{0}
\|_{L^{2}} e^{\|\mu\|_{L^\infty}\|u\|_{L^{1}(0,T)}}.
$$
This gives the continuity of the weak solutions with respect to the
initial conditions.

Now, let us assume that $\Psi_{0} \in H^{2}(\RR^N,\CC)$. Take $C$ to
be the bound of the multiplication operator $\mu$ over $H^2$: i.e.
$C$ is a positive constant such that for every $\varphi \in
H^{2}(\RR^N,\CC)$, $\|\mu\varphi\|_{H^{2}}$ $\leqslant C
\|\varphi\|_{H^{2}}$. We consider, then, $T_2>0$ such that  $C
\|u\|_{L^{1}(0,T_2)}$ $< 1$. By applying the fixed point theorem on
$$
\Theta_{2}:C^{0}([0,T_2],H^{2}) \rightarrow C^{0}([0,T_2],H^{2})
$$
defined by the same expression as $\Theta$, and using the uniqueness
of the fixed point of $\Theta$, we get that the weak solution is a
strong solution. The continuity with respect to the initial
condition of the strong solution can also be proved applying the
same arguments as in above.
\\

Finally, let us justify that the weak solutions take values in
$\SSS$. For $\Psi_{0} \in H^{2}$, the solution belongs to
$C^{1}([0,T],L^{2}) \cap C^{0}([0,T],H^{2})$ thus, the following
computations are justified
$$
\frac{d}{dt} \|\Psi(t)\|_{L^{2}}^{2} = 2 \Re \her{\frac{\partial
\Psi}{\partial t}}{\Psi} =0.
$$
Thus $\Psi(t) \in \SSS$ for every $t\in [0,T]$.

For $\Psi_{0} \in \SSS$, we get the same conclusion thanks to a
density argument and the continuity for the
$C^{0}([0,T],L^{2})$-topology of the weak solutions with respect to
the initial condition. $\square$

\textbf{Acknowledgments :} The author thanks K. Beauchard and J-M.
Coron for many helpful discussions.

\bibliographystyle{plain}

\end{document}